\documentclass[12pt]{article}
\usepackage{theorem} \usepackage{amsmath} \usepackage{amssymb} \usepackage{enumerate}
\usepackage[dvips]{graphicx}
\theoremstyle{plain}
   \newtheorem{Theorem}{Theorem} \newtheorem{Lemma}{Lemma}
   \newtheorem{Proposition}{Proposition}   \newtheorem{Corollary}{Corollary}
\theorembodyfont{\rmfamily}
   \newtheorem{Definition}{Definition}  \newtheorem{Example}{Example} 
   \newtheorem{Remark}{{\it Remark}}
 \def\endproof{\hfill$\square$\vspace{4mm}}

\def\a{\alpha} \def\calA{{\mathcal A}}
\def\b{\beta} \def\calB{{\mathcal B}} \def\Bur{\mathrm{Bur}}
  \def\d{\delta} \def\D{\Delta}
\def\e{\epsilon}\def\f{\varphi} \def\fix{\mathrm{Fix}}
\def\g{\gamma} \def\G{\Gamma}\def\Homeo{\mathrm{Homeo}}
 \def\id{\mathrm{id}}  \def\ind{\mathrm{ind}}
\def\int{\mathrm{Int\,}} \def\Iso{\mathrm{Iso}} \def\J{\bar{J}}
\def\calJ{{\mathcal J}}\def\calK{{\mathcal K}}
\def\l{\lambda}\def\L{\Lambda}\def\calL{{\mathcal L}} \def\LCM{\mathrm{LCM}}
 \def\m{\mu} \def\N{{\mathbb N}}\def\n{\nu} 
\def\o{\omega} \def\P{\Phi}  \def\calP{\mathcal{P}}\def\PM{\mathrm{PM}}
\def\R{{\mathbb R}}\def\r{\rho}\def\calR{{\mathcal R}}
\def\s{\sigma} \def\calS{{\mathcal S}}  \def\calT{\mathcal{T}}  
 \def\t{\tau} \def\tr{\mathrm{tr}\,}  \def\Z{{\mathbb Z}}

\begin{document}

\title{The generalized Lefschetz number of homeomorphisms on punctured disks}

\author{Takashi Matsuoka\\Department of Mathematics, Naruto University of Education\\
Naruto, Tokushima 772-8502 Japan}
\date{}
\maketitle
\begin{abstract}
We compute the generalized Lefschetz number 
of orientation-preserving self-homeomorphisms of a compact punctured disk,
using the fact that homotopy classes of these homeomorphisms
can be identified  with braids.
This result is  applied to study Nielsen-Thurston 
 canonical homeomorphisms on a punctured disk.
 We determine, for a certain class of braids,  the rotation number of the corresponding canonical
homeomorphisms on the outer  boundary circle.
As a consequence of this result on the rotation number, 
it is shown that the canonical homeomorphisms corresponding to some braids 
are pseudo-Anosov with associated foliations having no interior singularities.
\end{abstract}

\section{Introduction}
The generalized Lefschetz number is one of the central notions in Nielsen fixed point theory. 
The classical Lefschetz number $L(f)$ is a well-known homotopy invariant for
proving the existence of a fixed point of a continuous self-map $f:X\to X$ on a 
connected, finite cell complex $X$.
It coincides with the fixed point index of the whole set $\fix(f)$ of fixed points, and
hence the non-vanishing of this number implies that $f$ has a fixed point. 

The generalized Lefschetz number $\calL(f)$ is a refinement of the Lefschetz number
 obtained by decomposing the fixed point set $\fix(f)$
into finitely many equivalence classes called fixed point classes.
On the fundamental group $\pi_1(X)$, an equivalence relation, called the Reidemeister equivalence,
is defined  using the induced action $f_{\pi}$ of $f$. 
An equivalence class under this relation is called a Reidemeister class.
Then, a Reidemeister class is assigned to each fixed point, and
the set of fixed points to which a given Reidemeister class $\a$ is assigned is called the fixed point class
determined by $\a$.
The compactness of $X$ implies that 
there are only finitely many Reidemeister classes determining non-empty fixed  point classes. 
This fact enables us to define the generalized Lefschetz number $\calL(f)$ 
as the formal sum of the Reidemeister classes with each class being 
indexed by fixed point index of the corresponding fixed point class.
Hence,  $\calL(f)$ is not an integer, but an element of the free abelian group
$\Z \calR(f_{\pi})$ generated by the  set $\calR(f_{\pi})$ of Reidemeister classes.
The non-vanishing of the coefficient of a Reidemeister class $\a$ in $\calL(f)$ implies the existence
of a fixed point with $\a$ assigned. 
Thus, by computing the generalized Lefschetz number,
we can prove the existence of a fixed point corresponding to each term in $\calL(f)$.
The generalized Lefschetz number is a homotopy invariant, and 
the classical one $L(f)$ is obtained from $\calL(f)$ by summing up the coefficients.
See e.g. \cite{bro,jia83,jia05} for general references of Nielsen fixed point theory.  

Practically, the generalized Lefschetz number is useful  in studying fixed points  only in the case where it is computable.
Unfortunately, it is very difficult to compute it from the definition in general. 
The Reidemeister trace formula \cite{rei, wec, hus} provides a method to compute it.
The classical Lefschetz number $L(f)$ is known to satisfy the following trace formula: 
If $f$ is a cellular map, 
$L(f)$ is equal to the alternating sum of the traces of the action of $f$ on the chain groups of $X$.
Analogously, the generalized Lefschetz number $\calL(f)$ satisfies the Reidemeister trace formula:
$\calL(f)$ is equal to  the alternating sum of the Reidemeister traces, which are the traces of the action 
 of a lift $\tilde{f}$ on the chain groups of the universal cover of $X$.
Despite the existence of this formula, however, it is still difficult to make a detailed computation,
particularly in the case of fundamental group being infinite and non-abelian. 
In this case, the author does not know any example of concrete computations carried out on large classes of maps.

In this paper, we compute the generalized Lefschetz number 
for orientation-preserving self-homeomorphisms $f$ of a compact punctured disk
which preserve the outer boundary circle (Theorem 1).
Such homeomorphisms are of great importance in the topological study of 2-dimensional dynamical systems,
for they include the homeomorphisms which are obtained from orientation-preserving
disk homeomorphisms by the blow-up construction at a finite, interior invariant set 
(see e.g. \cite[Section 1.6]{boy}).
 We should note that our computation is not complete in the sense that
the problem of distinguishing Reidemeister classes is left unsettled.
This means that we shall obtain an element in the group ring $\Z \pi_1(X)$ which is mapped to
the generalized Lefschetz number under the projection from $\Z \pi_1(X)$ to $\Z \calR(f_{\pi})$.
Thus, our result  may be thought of as giving an ^^ ^^ upper bound'' of the generalized Lefschetz number.
Our computation utilizes the fact that the homotopy class 
(or equivalently the isotopy class)  of $f$ can be identified with a braid.
We show that a braid is designated by a finite sequence of positive integers, and  we shall 
compute the generalized Lefschetz number directly from this sequence.
 For surfaces with boundary, Fadell and Husseini showed in \cite{fh83} 
that the computation of the Reidemeister trace is reduced to  that in
the Fox free differential calculus on free groups.
Our result is obtained by carrying out this computation.
In \cite{mat93}, the author computed the image of the generalized Lefschetz number $\calL(f)$
under the projection
from $\Z \calR(f_{\pi})$ to the ring  $\Z[t,t^{-1}]$ of integer polynomials in the variable $t$ and its inverse.
The present result improves the computation there.  

It is a natural question whether our method is applicable 
in the general case  where $f$ may not preserve the outer boundary circle.
In this case,  $f$ is thought of as an orientation-preserving homeomorphism on a punctured sphere, and its homotopy class is identified with a braid on a sphere. 
Our method is based on the fact that
 a braid on a plane is designated by a finite sequence of positive integers. 
  At this moment, the author does not know a similar fact on a sphere,
and can not give an answer to the question.

On surfaces with boundary, Wagner \cite{wag} exploited an algorithm
 to compute the generalized Lefschetz number 
for a continuous map whose action on the fundamental group satisfies an algebraic condition.
This condition is satisfied by most of continuous maps, but not by homeomorphisms.
Therefore, the Wagner's algorithm is not applicable to our case.
 
We give two applications of our result in Section 4.
The Nielsen-Thurston classification theory of isotopy classes of surface homeomorphisms  
provides a representative $\f$, called a canonical homeomorphism in each isotopy class of surface homeomorphisms. 
Canonical homeomorphisms play an essential role
in the study of dynamics of surface homeomorphisms, because it
 has the ^^ ^^ simplest'' dynamical complexity among the homeomorphisms in its isotopy class.
For instance, all the periodic points of $\f$ persist under homotopy.
We apply our result on $\calL(f)$ to study periodic points of
canonical homeomorphisms on a punctured disk: We determine, for a certain class of braids,  
the rotation number of the corresponding canonical
homeomorphisms on the outer  boundary circle (Proposition 2).

The second application concerns the problem of determining the type of the
 canonical homeomorphism in a given isotopy class.
There is an algorithm to solve this problem due to Bestvina and Handel \cite{bh}.
 Similar algorithms for the disk case were given in \cite{fm,los}.
Also, different algorithms were given in \cite{bgn,hc}.
Our theorem provides an algebraic approach to this problem.
We show that our result on the rotation number on the outer boundary circle 
implies that the canonical homeomorphisms corresponding to some families of braids 
are pseudo-Anosov with associated foliations having no interior singularities (Proposition 3).

In the last section, as a by-product of an argument in the proof of Proposition 2,
we give a lower and an upper bound for the Nielsen number $N(f)$
 for the class of braids treated in Proposition 2. 

{\it Acknowledgement}. The author would like to  thank the referee for making helpful suggestions that improved exposition of the paper.
\section{Generalized Lefschetz number}

We recall  the definition of the generalized Lefschetz number. 
 Let $X$ be a connected finite cell complex, and $f:X\to X$ a continuous map.
Let $\fix(f)$ be the fixed point set of $f$.
Choose a base point $x_0$ of $X$, and let $\pi$ denote the fundamental group $\pi_1(X,x_0)$  of $X$.

 Given a homomorphism $\psi:\pi\to\pi$, 
two elements $\l_1,\l_2$ of $\pi $ are said to be {\it Reidemeister equivalent} with respect to $\psi$
(or {\it $\psi$-conjugate})
if there is a $\l\in \pi$ such that 
\[\l_2=\psi(\l)\l_1{\l}^{-1}.\]
An  equivalence class under this equivalence relation is called 
a {\it Reidemeister class}.
Let $\calR(\psi)$ denote the set of Reidemeister classes, and $\Z \calR(\psi)$  
the free abelian  group generated by the elements of $\calR(\psi)$. 

Choose a path $\t$ from  $x_0$ to $f(x_0)$. 
This is called a {\it base path}.
Let $f_{\pi}:\pi  \to \pi$  denote the composition 
of $f_*:\pi_1(X,x_0)\to \pi_1(X,f(x_0))$ with the isomorphism $\t^{-1}_*:\pi_1(X,f(x_0)) \to \pi_1(X,x_0)$
 induced by $\t^{-1}$.
We shall consider Reidemeister classes with respect to $f_{\pi}$.
 For $x \in  \fix(f)$, take  a path $l$ from the base point $x_0$ to $x$. Then,
it is easy to see that the Reidemeister class 
represented by $[\t(f\circ l)l^{-1}] \in \pi$ is independent of the choice of  $l$. 
This class is denoted by  $R(x)$ and is called the Reidemeister class (or a coordinate) of $x$.
For a Reidemeister class $\a\in \calR(f_{\pi})$, 
let $\fix_\a(f)=\{x\in\fix(f)\;| \; R(x)=\a\}$. 
This set is called the {\it fixed point class} of $f$ determined by $\a$.
We then have the decomposition 
\[\fix(f)=\bigcup_{\a\in\calR(f_{\pi})}\fix_{\a}(f).\]
The compactness of $X$ implies that $\fix_\a(f)$ is empty except for finitely many $\a$.
For an isolated set $S$ of fixed points of $f$, let $\ind(S)$ denote the fixed point index of $S$ 
with respect to $f$.
\begin{Definition}
The generalized Lefschetz number $\calL(f)$ of $f$ is defined by 
\begin{equation*}
\calL(f)=\sum_{\a\in \calR(f_{\pi})}\ind(\fix_\a(f))\,\a\; \in \;\;\Z\calR(f_{\pi}).
\end{equation*}
\end{Definition}@

The generalized Lefschetz number is a homotopy invariant in the following sense:
Let $g:X\to X$ be a continuous map  homotopic to  $f$ through a homotopy $\{h_t\}_{0\le t\le 1}$.
As a base path for $g$, take the composite of $\t$ with the path $h_t(x_0)$ $(0\le t\le 1)$ 
so that we have $f_{\pi}=g_{\pi}$. Then, the  Nielsen fixed point theory  asserts the equality
$\calL(f)=\calL(g)$. 

Let $\tilde{X}$ be the universal covering space of $X$.
For integers $q$, let $C_q(\tilde{X})$ be the $q$-chain group of $\tilde{X}$.
The action of $\pi$ on $\tilde{X}$ induces an action of the group ring $\Z \pi$ on $C_q(\tilde{X})$.
Then, $C_q(\tilde{X})$ becomes a finitely generated free $\Z \pi$-module.
If $f$ is a cellular map, its lift $\tilde{f}$ induces the twisted-module
homomorphism $\tilde{f}_{\sharp q}:C_q(\tilde{X})\to C_q(\tilde{X})$.
Then,  a trace $\tr \tilde{f}_{\sharp q}$ is defined as an element of $\Z \calR(f_{\pi})$.
The Reidemeister trace formula \cite{rei, wec, hus} asserts that
\[\calL(f)=\sum_{q\ge 0}(-1)^{q} \tr \tilde{f}_{\sharp q}.\]

 Note that the classical Lefschetz number is equal to the sum of the coefficients in $\calL(f)$, and  
the Nielsen number $N(f)$ is the number of Reidemeister classes with non-zero coefficients in $\calL(f)$.

\section{Computation on punctured disks}
We shall fix an integer $n$ with $n\ge 3$.
Let $D_n$ be a compact $n$-punctured disk, namely, it is a compact surface obtained from
 a closed disk $D$ by removing the interiors of $n$ disjoint closed disks $D(1),\ldots,D(n)$ contained in
the interior of  $D$.
$D_n$ has $n+1$ boundary circles.
One of these is  $\partial D$ called the {\it outer boundary circle} of $D_n$,
and the others $\partial D(1),\ldots,\partial D(n)$  are called the inner boundary circles of $D_n$. 
Let $\Homeo_+(D_n,\partial D)$ denote the set of orientation-preserving homeomorphisms 
$f:D_n\to D_n$ which preserve the outer boundary circle $\partial D$ setwise.
In this paper, we shall compute the generalized Lefschetz number $\calL(f)$ for
any $f\in \Homeo_+(D_n,\partial D)$ up to distinguishing Reidemeister classes.

An isotopy class of such homeomorphisms can be identified with a braid:
Let $\Iso_+(D_n,\partial D)$ be the group of isotopy classes of homeomorphisms in $\Homeo_+(D_n,\partial D)$.
Let $B_n$ denote the $n$-braid group.
Then,  we have an isomorphism of groups 
\begin{equation*}\Iso_+(D_n,\partial D)\to B_n/Z_n, \end{equation*}
where $Z_n$ is the center of $B_n$.
This identification is defined in the following way:
Choose an isotopy $\{f_t:D\to D\}_{0\le t\le 1}$ such that  $f_0=\id$ and that
$f_1$ coincides with $f$ on $D_n$.
The existence of such an isotopy is guaranteed using the well-known Alexander's trick. 
Then, the  subset $\bigcup_{0\le t\le 1} (f_t(D(1)\cup\cdots\cup D(n))\times \{t\})$
 of $D\times [0,1]$ consists of disjoint $n$ tubes.
 If we regard each tube as a string, we obtain an $n$-braid.
 We denote this braid by $\b(f)$, and call it the {\it braid of} $f$.
The element of the quotient group $B_n/Z_n$ represented by
 $\b(f)$  does not depend on the choice of an isotopy $\{f_t\}$. Thus we obtain a map
$\Iso_+(D_n,\partial D)\to B_n/Z_n$. It is known that this map becomes  an isomorphism.

We can assume that the centers of the sub-disks $D(1),\ldots,D(n)$
lie on a line in that order, hence so do the initial points of the braid $\b(f)$. 
For $i=1,\dots,n-1$, we denote by $\s_i$ the $i$-th elementary  braid,
in which the $i$-th string overcrosses the $(i+1)$-th string once
 and all other strings go straight from the top to the bottom.
The braid group $B_n$ admits a presentation with generators
 $\s_1,\dots,\s_{n-1}$ and defining relations (see e.g. \cite{bir74}):
\begin{align*}
 \s_i\s_j&=\s_j\s_i && \text{if}\; |i-j|\ge 2, \;1\le i,j\le n-1,\notag \\
\s_i\s_{i+1}\s_i&=\s_{i+1}\s_i\s_{i+1} && 1\le i\le n-2. 
\end{align*}  

Define $\r\in B_n$ by  $\r=\s_{n-1}\cdots\s_2\s_1$.
 Let $\theta$ be the {\it full-twist} $n$-braid defined by $\theta =(\s_1\s_2\cdots\s_{n-1})^n$.
 $\theta$ is a generator of  the center $Z_n$. 
In particular, it commutes with every braid.
Note that $\r^n$ is equal to $\theta$, since $\r=\D(\s_1\cdots \s_{n-1})\D^{-1}$, where 
$\D$ is a half-twist braid $(\s_1\s_2\cdots\s_{n-1})\cdots (\s_1\s_2)\s_1$.

For  a positive integer $i$, let $\b(i)=\s_1^i\r\in B_n$.
Let $d$ be a positive integer.
Given a sequence $I=(i_1,\ldots,i_d)$ of positive integers, define an $n$-braid $\b(I)$
by
\[\b(I)=\b(i_1)\cdots\b(i_d)=\s_1^{i_1}\r\cdots\s_1^{i_d}\r.\]

The following  proposition has been proved in \cite{mat93}.
We give here a slightly simplified proof.
\begin{Proposition}
Any braid is conjugate to a braid of the form $\theta^{\m}\b(I)$, where $\m$ is an integer and
$I$ is a finite sequence of positive integers.
\end{Proposition}

{\it Proof}.
By the defining relations of $B_n$, it is easy to see that for $i=1,\ldots,n-2$, 
\[\s_i\r=\r\s_{i+1}.\]
This implies that  
 \begin{equation}\s_i=\r^{1-i}\s_1\r^{i-1}\end{equation}
for any $ i$. Also, we have  
\begin{equation}(\s_1\r)^{n-1}=\theta. \end{equation}
For any  $i$, we have by (1), (2) 
\[\s_i^{-1}=\r^{1-i}\s_1^{-1}\r^{i-1}=\r^{2-i}(\s_1\r)^{-1}\r^{i-1}=\theta^{-1}\r^{2-i}(\s_1\r)^{n-2}\r^{i-1}.\]
This and (1) imply that $\s_i$'s and $\s_i^{-1}$'s 
 can be written as a product of $\s_1$, $\r$, $\r^{-1}$, and $\theta^{-1}$, and hence
any braid $\b$ is conjugate to a braid of the form $\theta^{\m}\s_1^{k_1}\r^{l_1}\cdots\s_1^{k_s}\r^{l_s}$,
where $\m\le 0$, $k_1,\ldots,k_s>0$ and $l_1,\ldots,l_s\in \Z$.
We can rewrite it in the form where all the exponents of $\r$ are equal to 1. In fact, 
since $\r^{-1}=\theta^{-1}(\s_1\r)^{n-2}\s_1$ by (2), we have 
$\r^j=\r^{nk}\r^{-l}=\theta^{k-l}((\s_1\r)^{n-2}\s_1)^l$ for any integer $j$,
where $k$ is an integer and $0\le l < n$ with $j=kn-l$.
\endproof

Note that the arguments in the proof also give a procedure how to find $\m, I$ and $\g$
 with $\b=\g^{-1}\theta^{\m}\b(I)\g$ for a given $\b\in B_n$.
\begin{Example}
Let $n=3$ and consider $\s_1\s_2^{-1}$. Since $\s_2^{-1}=\theta^{-1}(\s_1\r)\r$, we have 
$\s_1\s_2^{-1}=\theta^{-1}\s_1^2\r^2$.
$\r^2$ is equal to $(\s_1\r)\s_1$, since  $kn-l=2$ for $k=l=1$. 
Therefore,  $\s_1\s_2^{-1}=\theta^{-1}\s_1^3\r\s_1=\s_1^{-1}\theta^{-1}\b(4)\s_1$. 
Hence, $\m=-1, I=(4)$, and $\g=\s_1$.
\end{Example}
\begin{Remark}
 $\m$ and  $I$ in this proposition are not unique.
For instance, we have by (2)
\[\b(i,\underbrace{1,\ldots,1}_{n-2},j)=\s_1^{i-1}(\s_1\rho)^{n-1}\s_1^j\rho=\theta \b(i+j-1).\]
Also, in $B_3$, since (2) implies $\s_1\rho\s_1=\theta\rho^{-1}$ and hence
$(\s_1\rho\s_1)^2=\theta^2\r^{-2}=\theta\r$,
we have for $i,j\ge 2$
 \[\b(i,2,j)=\s_1^{i-1}(\s_1\rho\s_1)^2\s_1^{j-1}\rho=\theta\b(i-1,j-1).\]
\end{Remark}

The purpose of this paper is to compute the generalized Lefschetz number $\calL(f)$
in terms of $\m$, $I$, and $\g$ given in Proposition 1.
Let $d$ be a positive integer, and $\Z_{d}$ the set $\{1,\ldots,d\}$ of integers mod $d$.
To state our main result, it is necessary to introduce the notion of a partition of $\Z_d$.
\begin{Definition}  
\begin{enumerate}
\item For integers $1\le p,q\le d$,
define a sequence $[p,q]$ of consecutive integers mod $d$ by 
\[[p,q]=\begin{cases}(p,\ldots,q) & \text{if}\quad p\le q\\
(p,\ldots,d,1,\ldots,q) & \text{if} \quad p>q.\end{cases}\]
This sequence is called a {\it block} in $\Z_d$, and the number of integers  contained in it
is called its length. For a block $B$, let $\underline{B}$ denote its underlying set, 
the set of integers contained in $B$.

\item A set $\{B_1,\ldots,B_s\}$ of  blocks in $\Z_d$ is a {\it partition} of $\Z_d$ if 
\begin{enumerate}
\item the underlying sets  $\underline{B_1}, \ldots, \underline{B_s}$ are mutually disjoint 
and $\underline{B_1}\cup\cdots \cup \underline{B_s}=\Z_d$, and  

\item each of $B_1,\ldots,B_s$ has length less than or equal to $n-1$.
\end{enumerate}
Note that a partition contains at most one block of type $[p,q]$ with $p>q$.

\item  Let $\calP(d)$ denote the set of partitions of $\Z_d$.
\end{enumerate}
\end{Definition}

\begin{Example} Assume $n\ge 5$. Consider the case of $d=4$. In this case, any block 
has length less than or equal to $n-1$. Therefore, $\calP(4)$ consists of the following fifteen partitions:
\[\{(1),(2),(3),(4)\}\]
\[\{(1,2),(3),(4)\},\;\;\{(1),(2,3),(4)\},\;\;\{(1),(2),(3,4)\},\;\;\{(2),(3),(4,1)\},\]
\[\{(1,2),(3,4)\},\;\;\{(2,3),(4,1)\},\]
\[\{ (1,2,3),(4)\},\;\;\{(1),(2,3,4)\},\;\;\{(2),(3,4,1)\},\;\;\{(3),(4,1,2)\},\]
\[\{(1,2,3,4)\},\;\;\{(2,3,4,1)\},\;\;\{(3,4,1,2)\},\;\;\{(4,1,2,3)\}.\]
 \end{Example}

The fundamental group $\pi=\pi_1(D_n,x_0)$ is identified with  a free group $F_n$ of rank $n$.
We shall define an action of $B_n$ on $F_n$.
Let $e$ be the unit element of $F_n$. Assume that $x_0\in \partial D$.
Let $\xi_1,\ldots,\xi_n$ be the standard generators of $\pi=F_n$
which are defined in the following way: 
We can assume that $D$ is the disk in the plane $\R^2$ with center $(0,0)$ and radius $2$, $x_0=(0,2)$,
and for $i=1,\ldots,n$,
 the sub-disk $D(i)$ has radius $\frac{1}{2(n+1)}$ and center $(-1+\frac{2}{n+1}\, i, \,0)$.
Then, the element $\xi_i$ is represented by a loop which  traces a straight line connecting $x_0$ 
to a point in $\partial D(i)$, encircles $\partial D(i)$  once in the anti-clockwise direction,
and retraces the line back to $x_0$. 
An action of the braid group $B_n$ on $F_n$ is defined by putting 
 $\s_i(\xi_j)=\xi_i\xi_{i+1}\xi_i^{-1}$, $\xi_i$, or $\xi_j$ according to whether
$j=i,j=i+1$, or $j\ne i, i+1$.
Thus, any braid $\b$ can be thought of as an automorphism  of $F_n$ (see \cite[Corollary 1.8.3]{bir74}).

In  the case of $\b(f)$, the  corresponding automorphism of $F_n$  can be described geometrically by using the isotopy $\{f_t\}$. 
Let $Z=\bigcup_{0\le t \le 1}(f_t(D_n)\times \{t\})$.
$Z$ is a solid cylinder with  $n$ disjoint open tubes removed.
Define a vertical path $v$ in  $Z$ by $v(t)=(x_0,t)$.
For $\e=0,1$, let $i_{\e}:D_n\to Z$ be the inclusion map defined by $i_{\e}(x)=(x,\e)$.
Given an element  $w\in \pi$,  choose a loop $l$  based at $x_0$ representing $w$. 
Then, it is easy to see that the loop $v^{-1}(i_0 \circ l)v$ in $Z$ is homotopic to 
 $i_1\circ l'$ in $Z$ for some loop    $l'$ in $D_n$.
Then  the image $\b(f)(w)$ coincides with the element of $\pi$ represented by $l'$.

As a base path $\t$, we shall take $\t$ given by  $\t(t)=f_t(x_0)$. Then, for any loop $l$ in $D_n$ based at $x_0$, 
the loop  $v^{-1}(i_0\circ l)v$ is homotopic to $i_1\circ (\t (f\circ l) \t^{-1})$ in $Z$. 
Therefore, $\b(f)(w)\in \pi$ is represented by the loop $\t (f\circ l) \t^{-1}$, 
and hence it is equal to $f_{\pi}(w)$.
 Thus, we have shown that
 \begin{equation}f_{\pi}=\b(f):F_n\to F_n.\end{equation}
For $w\in F_n$, we shall use the symbol $w^\b$ to denote its image under the automorphism $\b$.

In our computation, we shall not use the standard generators, but use the 
 generators $a_1,\dots,a_n$ for $F_n$ defined  by $a_i=\xi_1\cdots \xi_i$.
Then, the action of $\s_i$  on  $F_n$ is written in a slightly simpler way as 
\begin{equation}
a_j^{\s_i}=\begin{cases}
       a_{i+1}a_i^{-1}a_{i-1}  &  \text{if}\;\; j=i, \notag \\
     a_j &   \text{if}\;\; j\ne i, \end{cases}
\end{equation}
where we put $a_0=e$. 
Note that $a_n^{\b}=a_n$ for any braid $\b$, since $a_n^{\s_i}=a_n$ for any $i$.

Let $\Z F_n$ be the group ring of $F_n$ over $\Z$. 
For $\b\in B_n$, the automorphism $\b$ of $F_n$  induces the ring automorphism of $\Z F_n$, 
which will be denoted by the same letter $\b$.
For $\eta\in \Z F_n$ and $\b\in B_n$, let $\eta^{\b}\in \Z F_n$ denote the image of $\eta$ under $\b$.

Let  $I=(i_1,\ldots,i_d)$ be a sequence of 
positive integers.
We shall introduce a map $W_I:\calP(d)\to \Z F_n$ which is necessary to state the main result.
First, for integers $j\ge 0$Cdefine $c_j\in F_n$ and $g_j\in \Z F_n$ by
\begin{equation}c_j=\begin{cases}a_2^{j/2} & \text{if\;} j \;\text{is even},\notag \\
a_1a_2^{(j-1)/2} & \text{if}\;j \;\text{is odd,}\end{cases}
 \end{equation}
 \[g_j=(-1)^{j+1}c_j.\]
 For $1\le l\le d$, let $\b_l(I)=\b(i_l,\ldots,i_d)\in B_n$. Note that $\b_1(I)=\b(I)$.
Suppose a block $B=[p,q]$  in $\Z_d$ is given.
Denote its length  by  $|B|$.
Define the braids $\a(B),\o(B)$ by
\begin{equation}\a(B)=\b_p(I),\quad \o(B)=\begin{cases} \b_q(I) &  \text{if}\quad  p\le q, \notag \\ 
 \b_q(I)\b(I)^{-1} &  \text{if}\quad p>q,
 \notag\end{cases}\end{equation}
and define $W_I(B)\in \Z F_n$ as follows:
\begin{equation}
W_I(B)=\begin{cases}
(g_0+\cdots+g_{i_p-2})^{\a(B)}a_{|B|+1}^{\o(B)} & \text{if}\quad |B|< n-1,\;i_p\ge 2, \notag \\
0 &  \text{if}\quad |B|< n-1,\; i_p=1,\notag\\
g_{i_p}^{\a(B)}a_{n-1}^{\o(B)}& \text{if} \quad |B|=n-1.
\end{cases}
\end{equation}
Then, the map $W_I:\calP(d)\to \Z F_n$ is defined  as follows:
Let  $\calB\in \calP(d)$. 
Let $B_1=[p_1,q_1],\ldots, B_s=[p_s,q_s]$ be the blocks in  $\calB$.
We can assume that $1\le p_1<p_2<\cdots <p_s\le d$ by rearranging the blocks if necessary.
Then, define $W_I(\calB)$ by
\[W_I(\calB)=W_I(B_1)\cdots W_I(B_s).\]

 Let  $\P_{\b}:\Z F_n\to \Z \calR(\b)$ denote the surjective homomorphism  induced by the projection
 $F_n=\pi\to \calR(\b)$.
By the definition of the Reidemeister equivalence, we have
\begin{equation}\P_{\b}(w)=\P_{\b}(w^{\b}) \quad \text{ for any}\; w\in F_n.\end{equation}
More generally, we have 
\begin{equation}\P_{\b}(w'w)=\P_{\b}(w^{\b}w') \quad \text{ for any}\; w, w'\in F_n.\end{equation}

Recall that $\b(f)$ can be written as $\g^{-1}\theta^{\m}\b(I)\g$ for some $\m\in\Z, \g\in B_n$,
and some sequence $I$ of positive integers.
Our main result is the following:
\begin{Theorem} Suppose $f:D_n\to D_n$ is an orientation-preserving homeomorphism 
which preserves the outer boundary circle setwise.
We choose an isotopy
$\{f_t\}:D\to D$ such that $f_0=\id$ and  $f_1$ coincides with $f$ on $D_n$.
As a base path  for $f$, take the path $\t$ defined by $\t(t)=f_t(x_0)$.
  Suppose  $\b(f)=\g^{-1}\theta^{\m} \b(I)\g$, 
where $\m$ is an integer, $\g\in B_n$, and  $I$ is a sequence of
positive integers with length $d$. Then 
\[\calL(f)=-\P_{\b(f)}\biggl(a_n^{\m}\sum_{\calB\in \calP(d)}W_I(\calB)^{\g}\biggr)\;\in \Z \calR(\b(f)).\]
\end{Theorem}

\begin{Example}
(a) Let $\b(f)=\b(i)$, where $i\ge 2$.
In this case,  $\m=0$, $\g=e$, $I=(i)$, and $d=1$. 
The partition $\{(1)\}$ is the only element of $\calP(d)=\calP(1)$,
 and    $\a((1))=\o((1))=\b(i)$. 
Therefore, by the above theorem and (4), we have 
\begin{align*}\calL(f)&=-\P_{\b(f)}(W_{(i)}(\{(1)\}))=-\P_{\b(i)}((g_0+\cdots +g_{i-2})^{\b(i)}a_{2}^{\b(i)})\\
&=-\P_{\b(i)}(g_2+\cdots + g_i)\\
&=\P_{\b(i)}(c_2)-\P_{\b(i)}(c_3)+\cdots +(-1)^{i}\P_{\b(i)}(c_i).\end{align*}

(b) Let $\b(f)=\b(i_1,i_2)$, where $i_1,i_2\ge 2$.
In this case, $\m=0$, $\g=e$, $I=(i_1,i_2)$, and $d=2$.
 $\calP(d)=\calP(2)$ consists of the three partitions
$\calB_1=\{(1),(2)\}, \calB_2=\{(1,2)\}$, and $\calB_3=\{(2,1)\}$.
Therefore, we have $\calL(f)=-\P_{\b(f)}(W_I(\calB_1)+W_I(\calB_2)+W_I(\calB_3))$, where
\[W_I(\calB_1)=\sum_{j=2}^{i_1}\sum_{k=2}^{i_2} g_{j}^{\b(I)}g_k^{\b(i_2)},\]
\[W_I(\calB_2)=\begin{cases}\Bigl(\sum_{j=0}^{i_1-2}g_{j}^{\b(I)}\Bigr)a_{3}^{\b(i_2)}
& \text{if}\; n\ge 4,\\
g_{i_1}^{\b(I)}a_2^{\b(i_2)} & \text{if}\; n=3,
\end{cases}\]
\[W_I(\calB_3)= \begin{cases}
\Bigl(\sum_{k=0}^{i_2-2}g_{k}^{\b(i_2)}\Bigr)a_{3} & \text{if}\; n\ge 4,\\
g_{i_2}^{\b(i_2)}a_2 & \text{if}\; n=3.\end{cases}\]
 \end{Example}

As a consequence of our  theorem, we can give an upper bound for the Nielsen number $N(f)$.
For $\eta\in \Z F_n$, let $\n(\eta)$ denote the number of elements 
of $F_n$ with non-zero coefficient in $\eta$.
Then, for a block $B=[p,q]$, we have 
\[\n(W_I(B))=\begin{cases}
i_p-1 &\text{if}\quad |B|<n-1, \\
1& \text{if}\quad |B|=n-1.\end{cases}\]
For a partition  $\calB=\{B_1,\ldots,B_s\}$, let $\n_I(\calB)=\n(W_I(B_1))\cdots \n(W_I(B_s))$.
Then, we have 
\begin{Corollary} Under the same hypothesis of Theorem 1, we have 
$N(f)\le  \sum_{\calB\in \calP(d)}\n_I(\calB)$.
\end{Corollary}
{\it Proof}.
It follows from Theorem 1 that
\[N(f)=\n(\calL(f))\le \n\Bigl(\sum_{\calB\in \calP(d)}W_I(\calB)\Bigr)
\le \sum_{\calB\in \calP(d)}\n(W_I(\calB)).\]
Since $\n(W_I(\calB))\le \prod_{r=1}^s\n(W_I(B_r))=\n_I(\calB)$, this gives the proof. 
\endproof
\begin{Example}
Consider the braids treated in Example 3.

(a) Let $\b(f)=\b(i)$, $i\ge 2$. Then, the partition $\{(1)\}$
is the only element of $\calP(1)$,
and hence  $\sum_{\calB\in \calP(d)}\n_I(\calB)=\n_{(i)}(\{(1)\})=i-1$. 
Therefore, by Corollary 1, we have $N(f)\le i-1$.
In fact, the equality  $N(f)=i-1$ holds for this braid.
This is proved by using the ring homomorphism  $\calT:\Z F_n\to \Z[t,t^{-1}]$   defined by 
$\calT(\xi_j)=t$ for $j=1,\ldots,n$. Then, $\calT$ induces the homomorphism $\calT:\Z \calR(\b(f))\to \Z[t,t^{-1}]$,
and the Reidemeister classes  $\P_{\b(i)}(c_2),\ldots, \P_{\b(i)}(c_{i})$ in $\calL(f)$ are sent to 
 mutually different elements $t^2,\ldots,t^i$ by $\calT$.
 Hence  they are different elements in $\calR(\b(i))$ and we obtain $N(f)=i-1$.

(b) Let $\b(f)=\b(i_1,i_2)$, where $i_1,i_2\ge 2$. Consider first the case of $n\ge 4$.
Let $\calB_1, \calB_2, \calB_3$ be the partitions as in Example 3(b).
Then, 
 $\n_I(\calB_j)=(i_1-1)(i_2-1), i_1-1, i_2-1$ for $j=1,2,3$ respectively, and so
we have  $\sum_{\calB\in \calP(d)}\n_I(\calB)=(i_1-1)(i_2-1)+(i_1-1)+(i_2-1)=i_1i_2-1$.
Hence, it follows from Corollary 1 that $N(f)\le i_1i_2-1$.
Consider next the case of $n=3$. Then, Corollary 1 implies $N(f)\le (i_1-1)(i_2-1)+2$.
If $i_1,i_2\ge 3$, the sharper estimate $N(f) \le (i_1-1)(i_2-1)-2$ holds, because the images of 
$W_I(\calB_2)=-g_{i_1}^{\b(I)}g_2^{\b(i_2)}$ and
$W_I(\calB_3)=-g_{i_2}^{\b(i_2)}g_2$ under $\P_{\b(f)}$ cancel by the images  of two terms in $W_I(\calB_1)$.
For a class of braids including this example, we shall give a sharper estimation  than Corollary 1 in Section 8.
\end{Example}

\begin{Remark} The image of  $-\calL(f)$ under $\calT$ coincides with 
the trace of the reduced Burau matrix $\Bur(\b(f))$ of the braid $\b(f)$ (cf. \cite{hj}) .
 This trace was computed in \cite{mat93} using the same expression
of braids as in Proposition 1.
Given  a square matrix $A$ of size $\n$ with entries in a commutative ring $R$,
let $\PM(A;k)$ be the sum of principal minors of $A$ of order $k$ if $1\le k\le \n$ and zero otherwise.
Then, we have the equality
\[\tr A^d=\sum_{\calB\in\calP(d)}(-1)^{d+\sharp \calB}\PM(A;|B_1|)\cdots \PM(A;|B_s|)\]
for any positive integer $d$, where $\calB=\{B_1,\ldots, B_s\}$.
Applying this to the case of $A=\Bur (\b(i))$, we have  for $I=(i,\ldots,i)\in \N^d$ that
  \[\tr \Bur(\b(I))=\tr \,\Bur(\b(i))^d=\sum_{\calB\in\calP(d)}(-1)^{d+\sharp \calB}P(i;|B_1|)\cdots P(i;|B_s|),\]
where   $P(i;k)$ denotes $\PM(\Bur(\b(i));k)$ for any $k$.
Theorem 1 in \cite{mat93} generalizes this 
equality  to an arbitrary sequence $I\in \N^d$ as follows:
\[\tr \Bur(\b(I))=\sum_{\calB\in\calP(d)}(-1)^{d+\sharp \calB}P(i_{p_1};|B_1|)\cdots P(i_{p_s}; |B_s|),\]
where $p_r$ is the initial element of $B_r$ for $1\le r\le s$.
Our main result, Theorem 1 above, gives a refinement of this equality.
\end{Remark}

\begin{Remark} In our setting, Reidemeister classes can be visualized by the method of 
Jiang\cite{jia88} using the mapping torus.
For $t\in [0,1]$, let $[t]$ denote the corresponding point in the circle $S^1=\R/\Z$.
Define a subspace  $T$ of $D\times S^1$ by $T=\bigcup_{0\le t\le 1} (f_t(D_n)\times \{[t]\})$,
which is homeomorphic to  the mapping torus of $f$. 
Then the set of Reidemeister classes is in one-to-one correspondence with
the set of free homotopy classes of loops in $T$.
The Reidemeister class $R(x)$ of $x\in \fix(f)$ corresponds to the free homotopy class
of the loop $(f_t(x),[t])\;(0\le t\le 1)$ under this identification.   
\end{Remark}

\section{Nielsen-Thurston classification of surface homeomorphisms}
We shall apply the theorem in the previous section to study periodic points of Nielsen-Thurston 
 canonical homeomorphisms
on a punctured disk, and also to 
the classification problem of homeomorphisms into isotopy classes.
We recall briefly the Nielsen-Thurston classification theory of surface homeomorphisms(\cite{flp,th}).
Let $M$ be a compact surface.
 A homeomorphism $\f:M \to M$ is said to be of {\it finite order} if some  of its  iterates
is equal to the identity map. 
The map   $\f$ is said to be 
{\it pseudo-Anosov}, if the following conditions are satisfied:
\begin{enumerate}
\item There exists a pair of transverse foliations on $M$, 
carrying measures which are uniformly expanded and  contracted by $\f$ respectively.
\item These foliations have  finitely many singularities which coincide in the 
interior $\int M$ and alternate on the boundary $\partial M$.
Any singularity is $p$-pronged for some integer $p\ge 3$ if it is
in the interior of $M$, and it is 3-pronged if it is in $\partial M$.
(We consider segments of the boundary to be prongs.)
\end{enumerate}
$\f$  is said to be {\it reducible}  
if there exists a finite collection of disjoint
annuli in $M$ such that $\f$ maps  their union $A$ to itself, and that 
each connected component $N$ of $M-A$, called a {\it component} of $\f$,  has negative Euler characteristic and for
any iterate $\f^m$ mapping $N$  to itself,
 its restriction to $N$ is either of finite order or pseudo-Anosov. 
 The Nielsen-Thurston  classification theory states that
every homeomorphism $f:M\to M$   is isotopic 
to a  homeomorphism $\f:M\to M$ which is of finite order, pseudo-Anosov, or reducible. 
The homeomorphism $\f$ is called a  {\it canonical homeomorphism} in the isotopy class of $f$.
In the case where $\f$ is irreducible, i.e.,  of finite order or pseudo-Anosov, the surface 
$M$ is called the component of $\f$.

One of the  common features of canonical homeomorphisms  $\f$ is that they have periodic points
on every boundary circle $C$. 
In fact, in the case where  $C$ is contained in a pseudo-Anosov component, 
the singularities of associated foliations on it are periodic points. 
Also,  in the case where  $C$ is contained in a finite-order component $N$, 
all the points  in $N$ are   periodic.
Since the restriction of $\f$ to $C$ is an orientation-preserving homeomorphism
of a circle, the periodic points in $C$  have the same least period.
We shall consider the problem of determining the  period of periodic points and
 the rotation number on $C$ in the case  where  $M=D_n$  and $C=\partial D$.
The reason why we choose the outer boundary  circle as the subject is that 
 this is the easiest case  to deal with by using the generalized Lefschetz number.
The result we shall obtain  will be applied to classify homeomorphisms  up to isotopy.

 Let $\f$ be an orientation-preserving canonical homeomorphism on $D_n$ preserving 
$\partial D$ setwise.
We denote by  $m(\f)$ the least period of periodic points on $\partial D$.
Let $N_{\f}$ be the component of $\f$ containing $\partial D$.
Choose an isotopy $\f_t:D\to D$ such that  $\f_0=\id$ and that 
$\f_1$ coincides with $\f$ on ${D_n}$. 
Assume the base point $x_0$ is in $\partial D$.
Define a base path $\tau$ for $\f$ by $\t(t)=\f_t(x_0)$.
Note that $\tau$ is contained in $\partial D$.
For every positive integer $m$, define a base path $\t_m$ for $\f^m$ 
 by $\t_m=\t(\f\circ \t)\cdots(\f^{m-1}\circ \t)$. Choose  a periodic point $x$ on $\partial D$.
Since $x_0$ and $x$ are contained in $\partial D$, we can choose  a path $l$ connecting these points
contained in $\partial D$. Then, the loop 
$\t_{m(\f)}(\f^{m(\f)}\circ l)l^{-1}$ is contained in $\partial D$ 
and hence it represents an element $a_n^{\n(\f)}\in F_n$ for some integer $\n(\f)$.
Note that  $\n(\f)$ does not depend on the choice of the periodic point $x$ and the path $l$.
It depends, however, on the choice of an isotopy $\f_t$, but is uniquely determined modulo $m(\f)$.
The number $\n(\f)/m(\f)$ modulo $\Z$ is equal to the  rotation number of $\f$ on $\partial D$.

The following lemma shows that, in the case where $m(\f)$ and $\n(\f)$ are relatively prime,
   the problem of determining these numbers  is reduced to 
the computation of the generalized Lefschetz number.
\begin{Lemma} Let $m$ and $\n$ are integers with $m>0$.
Assume that they are relatively prime, and that
  the coefficient of $\P_{\b(\f^m)}(a_{n}^{\n})\in \calR(\b(\f^m))$ in $\calL(\f^{m})$ is non-zero. 
Then, $m=m(\f)$ and $\n=\n(\f)$.
\end{Lemma}

{\it Proof}.
Since  the coefficient of $\P_{\b(\f^m)}(a_{n}^{\n})$ in $\calL(\f^m)$ is non-zero,
$\f^{m}$ has a fixed point $x$ with Reidemeister class $\P_{\b(\f^m)}(a_n^{\n})$.
Then, we can take a path $l$ from the base point $x_0$ to $x$ so that
$\t_m(\f^m\circ l)l^{-1}$ represents $a_n^{\n}$. 
We shall show that $x$ is $\f^m$-related to $\partial D$, namely
there exists a  path connecting a point in $\partial D$ with  $x$ which is homotopic to its image under $\f^m$
via a homotopy of paths  such that each path in the homotopy  connects a point in $\partial D$ with  $x$.
Choose a loop $\l$ contained in $\partial D$ based at $x_0$ which represents $a_n^{\n}$.
Let $\eta=\t_m^{-1}\l$ and for $0\le u\le 1$, let $\eta_u(t)=\eta(ut+1-u)$ and $l_u=\eta_u l$. 
Then, $\{l_u\}$ gives a homotopy of  paths such that $l_u(0)\in \partial D, l_u(1)=x$, and
$l_0$ and $l_1=\t_m^{-1}\l l$ are homotopic to $l$ and $\f^m\circ l$ fixing end points respectively.
Thus, we have shown that $l$ is a desired path to prove $x$ being $\f^m$-related to $\partial D$.
(This is proved in a different way by  Guaschi \cite[Proposition 14(b)]{gua}.)
Then, it follows from Jiang and Guo \cite[Lemma 3.4]{jg} that
$x\in\partial D$ if $\f|_{N_{\f}}$ is  pseudo-Anosov,
 and there exists a path in $\fix(\f^m)\cap N_{\f}$ connecting $x$ to $\partial D$
 if $\f|_{N_{\f}}$ is of finite order.
Hence, $m=qm(\f)$ for some positive integer $q$.
Moreover, since $a_n$ is fixed under $\f_{\pi}=\b(\f)$, we have $\n=q\n(\f)$.
Since $m$ and $\n$ are relatively prime, $q$ must be one. Thus the proof is completed.
\endproof

Let $\LCM$ denote the least common multiple for positive integers.
Using this lemma, Theorem 1 can be applied to obtain 
\begin{Proposition} Suppose  the braid $\b(\f)$ is conjugate to $\theta^{\m}\b(I)$, where
$I=(i_1,\ldots,i_d)$ is a sequence of positive integers.
 Assume  either that $n \ge 4$ and $i_1,\ldots,i_d\ge 2$, or
that $n=3$ and $i_1,\ldots,i_d\ge 3$.
 Then \[m(\f)=\frac{\LCM\{d,n-2\}}{d}, \quad \n(\f)=m(\f)\m+\frac{\LCM\{d,n-2\}}{n-2}.\]
\end{Proposition}
This proposition  will be proved in Section 7 using some lemmas on the computation of the
generalized Lefschetz number given in Section 6. 
When $n=3$, this result cannot be extended to the case of $i_1,\ldots,i_d\ge 2$. 
For instance, $\b(2)\in B_3$ is conjugate to $\rho^2$, and so $\f$ is of  finite-order  and $m(\f)=3$.
This is not equal to  $\displaystyle{\frac{\LCM\{d,n-2\}}{d}}=1$.

This proposition has a consequence on  the classification problem of canonical homeomorphisms on a punctured disk.
Boyland \cite{boy84} proved that if $n$ is prime and $\b(\f)$ is cyclic, that is,
the permutation on the 
punctures induced by $\b(\f)$ is cyclic, then  $\f$ is irreducible.
He also proved that if $\f$ is irreducible, $\b(\f)$ is cyclic, and the exponent sum of $\b(\f)$ is
not divisible by $n-1$,  then $\f$ is pseudo-Anosov.
In particular, if $n$ is prime, $\b(I)$ is cyclic, and $i_1+\cdots +i_d$ is not divisible by $n-1$, then
$\f$ is pseudo-Anosov.
Matsuoka \cite{mat93} has proved that, under the assumption of Proposition 2,
 the canonical homeomorphism $\f$ with braid $\b(I)$
 contains a pseudo-Anosov component, except only for the case where $n\ge 4,I=(2,\ldots,2)$ and $n=3, I=(3,\ldots,3)$.
 This result was proved by using the computation of the reduction  $\calT(\calL(f))$ mentioned in Remark 2. 
Our main theorem on the computation of  the unreduced number $\calL(f)$ can be applied to improve
 this result. In fact, as a consequence of Proposition 2, we have the following proposition.

\begin{Proposition} Assume $n\ge 5$. Let $I$ be a sequence of integers 
 $i_1,\ldots,i_d\ge 2$ which  are all odd or all even.
Assume that $n-2$ and  $d$ are relatively prime.
Then, the canonical homeomorphism $\f$ with braid  $\b(I)$  is pseudo-Anosov. 
Moreover, the foliations associated to $\f$ have no interior singularities.
\end{Proposition}

{\it Proof}.
Since $n-2$ and $d$ are assumed to be relatively prime,  we have $m(\f)=n-2$ and $\n(\f)=d$ by Proposition 2.
This implies that  the periodic points on $\partial D$ have period $n-2$ and rotation number $d/(n-2)$. 
Let $\m$ be the permutation  on the inner boundary circles of $D_n$ induced by $\f$.
Assume $\f|_{N_{\f}}$ were of finite-order. Then, $\f|_{N_{\f}}$  is topologically conjugate to  the 
rigid rotation on the unit disk by angle $2\pi d /(n-2)$ restricted to the exterior of an appropriate
set of punctures. 
Hence, there exist $n-2$ boundary circles $C_1,\ldots,C_{n-2}$ of $N_{\f}$ cyclically permuted by $\f$.
If none of $C_1,\ldots,C_{n-2}$ is a  boundary circle of $D_n$, each of them
 surrounds at least two boundary circles of $D_n$.
Therefore, there must exist at least $2(n-2)$ boundary circles of $D_n$.
 Since $n\ge 5$, this number exceeds $n$, which is impossible.
Therefore,  some of $C_1,\ldots,C_{n-2}$ is a  boundary circle of  $D_n$,
and so are all of $C_1,\ldots,C_{n-2}$, 
since they are cyclically permuted by $\f$.
Therefore, $\m$ has  a cycle with length $n-2$.
We shall show that this contradicts to an assumption of the proposition.
In the case where $i_1,\ldots,i_d$ are all even, $\m$ is equal to the permutation 
induced by $\r^d$,  and hence it is the $d$-th power of a cyclic permutation on $n$ circles.
Hence, $n-2$ must divide  $n$, which is a contradiction since $n\ge 5$.
Also, in the case where $i_1,\ldots,i_d$ are all odd, $\m$  fixes  one of the inner boundary circles
and on the other $n-1$ inner boundary circles, $\m$ is the permutation induced by
$(\s_1\r)^d$, which is the  $d$-th power of a cyclic permutation.
Thus,  $n-2$ must divide $n-1$, which is a contradiction. Therefore, $\f|_{N_{\f}}$ is not of finite order, and hence
it must be pseudo-Anosov.

Let $c$ be the number of inner boundary circles of $N_{\f}$.
Choose one of the foliations on $N_{\f}$ and let $\calS$ denote the set of its singularities. 
 Denote by $p(x)$ the number of prongs at a singularity $x$.
Then we have the following Euler-Poincar\'e formula
(see e.g. \cite{flp}, p.75):
\begin{equation}\sum_{x\in \,\calS}(2-p(x)) =2\chi(N_{\f})=2(1-c).\end{equation}
Since the singularities on $\partial D$ are periodic points with least period $n-2$,
there exist  at least $n-2$ singularities on $\partial D$.
Also, each inner boundary circle of $N_{\f}$ contains at least one singularity.
Hence $\sharp (\calS \cap \partial N_{\f})\ge  n-2+c$.
Therefore, since $2-p(x)=-1$ for every singularity on $\partial N_{\f}$ and $n\ge c$, we have by (6)
\begin{align*}
\sum_{x\in \calS\cap \int N_{\f}}(2-p(x))& 
=\sum_{x\in\calS}(2-p(x))-\sum_{x\in \,\calS\cap \partial N_{\f}}(2-p(x))\\
&=2(1-c) -(-\sharp(\calS\cap \partial N_{\f}))\ge n-c\ge 0.
\end{align*}
This implies that  there are no interior singularities on $N_{\f}$, 
since $2-p(x)<0$ for any $x\in \calS \cap \int N_{\f}$, and also that $0\ge n-c$.
Hence $c=n$, and so $N_{\f}=D_n$ and $\f$ is pseudo-Anosov.
\endproof.

The above proposition cannot be extended to the case of $n=3,4$. 
For instance, $\b(2)\in B_3$ is conjugate to $\rho^2$, which corresponds to
a finite-order homeomorphism. Also, $\b(2)\in B_4$ is conjugate to  $\r\s_3\s_2$,
which corresponds to a reducible homeomorphism having only  finite-order components.

\section{Proof of Theorem 1}
For surfaces with boundary, Fadell and Husseini showed in \cite{fh83} 
that the computation of the generalized Lefschetz number is reduced to that in
the Fox free differential calculus on free groups.
The Fox partial derivative operator  ${\partial}/{\partial a_j}:\Z F_n\to \Z F_n$, $j=1,\dots,n$,
 is defined by the following rules (see \cite{bir74,mor}):

\begin{enumerate}
\item $\displaystyle{\frac{\partial}{\partial a_j}(\eta_1+\eta_2)=\frac{\partial \eta_1}{\partial a_j}+ 
\frac{\partial \eta_2}{\partial a_j}, \quad \eta_1,\eta_2\in \Z F_n}$,
\item $\displaystyle{\frac{\partial}{\partial a_j}(w_1w_2)
=\frac{\partial w_1}{\partial a_j}+w_1\frac{\partial w_2}{\partial a_j}, \quad w_1,w_2 \in  F_n}$,
\item $\displaystyle{\frac{\partial a_i}{\partial a_j}=\d_{i,j}, \quad 1\le i,j\le n}$, 
 $\quad$ 

where $\d_{i,j}=1$ or $0\,$ 
according to whether $\,i=j$ or $i\ne j$.
\item $\displaystyle{\frac{\partial e}{\partial a_j}=0}$.
\end{enumerate}
These rules imply that for $v,w\in F_n$,
\begin{equation} 
\frac{\partial}{\partial a_j}(vwv^{-1})=
(1-vwv^{-1})\frac{\partial v}{\partial a_j}+v\frac{\partial w}{\partial a_j}.
\end{equation}

Given a braid $\b\in B_n$, let $J(\b)$ be the Jacobian matrix $( {\partial a_i^{\b}}/{\partial a_j})$.
As an application of the Reidemeister trace formula, 
Fadell and Husseini  proved that
$\calL(f)=\P_{\b(f)}(1-\tr J(\b(f)))$ (\cite[Theorem 2.3]{fh83}).
For a matrix $A$ with entries  in $\Z F_n$, 
let $A^{\b}$ denote the matrix obtained from $A$ by replacing each entry 
with its image  under $\b$. Then, we have by (5)
\begin{equation}\P_{\b}(\tr A'A)=\P_{\b}(\tr A^{\b}A')
\end{equation} 
for any matrices $A,A'$.
Using the chain rule for the Fox calculus,
we have $J(\b\b')=J(\b)^{\b'}J(\b')$
 for any braids $\b, \b'$.
 Let $\b\in B_n$. Since $a_n^\b=a_n$, 
the last row of $J(\b)$ is $(0\cdots 01)$.
Let $\J(\b)$ denote the reduced matrix obtained from $J(\b)$ by deleting the last column and the last row.
Then, $\tr \J(\b)=\tr J(\b)-1$ and  hence we have
\begin{equation}\calL(f)=-\P_{\b(f)}(\tr \J(\b(f))).
\end{equation}

We shall show that $\calL(f)$ is  determined essentially by $\b(I)$ .
Note that \begin{equation}\J(\b\b')=\J(\b)^{\b'}\J(\b'). \end{equation}
Note also that since $a_i^{\theta}=a_na_i a_n^{-1}$ for any $i$, we have $\J(\theta)=a_n I_{n-1}$, 
where $I_{n-1}$ is the identity matrix. Therefore, we have  that
\begin{equation} \J(\theta^{\m}\b)=\J(\theta^{\m})^{\b}\J(\b)=a_n^{\m}\J(\b).
\end{equation}
We have by (10) that \[\J(\g\b(f))=\J(\theta^{\m}\b(I)\g)=\J(\theta^{\m}\b(I))^{\g}\J(\g).\]
Also, we have $\J(\g) \J(\g^{-1})^{\g}=\J(\g\g^{-1})^{\g}=I_{n-1}$.
Therefore, using (8), we have 
\[\begin{split}\P_{\b(f)}(\tr \J(\b(f)))&=\P_{\b(f)}(\tr \J(\g^{-1})^{\g\b(f)}\J(\g\b(f)))\\
&=\P_{\b(f)}(\tr \J(\g\b(f))\J(\g^{-1})^{\g})\\
&=\P_{\b(f)}(\tr\J(\theta^{\m}\b(I))^{\g}\J(\g)\J(\g^{-1})^{\g})\\
&=\P_{\b(f)}(\tr \J(\theta^{\m}\b(I))^{\g}).
\end{split}\]
Therefore, we have by (9) and (11) that 
\begin{equation}\calL(f)=-\P_{\b(f)}(a_n^{\m}\,\tr \J(\b(I))^{\g}).\end{equation} 

\begin{Lemma} 
\begin{enumerate}
\item
Two elements  $w_1, w_2$ of $F_n$ are Reidemeister equivalent with respect to $\b(I)$ 
if and only if   $a_n^{\m}w_1^{\g}, a_n^{\m}w_2^{\g}\in F_n$ are Reidemeister equivalent with respect to $\b(f)$.
 \item 
 Suppose  $\eta_1,\eta_2 \in \Z F_n$.
 Then,  $\P_{\b(I)}(\eta_1)=\P_{\b(I)}(\eta_2)$ if and only if 
$\P_{\b(f)}(a_n^{\m}\eta_1^{\g})=\P_{\b(f)}(a_n^{\m}\eta_2^{\g})$.
\end{enumerate}
\end{Lemma}

{\it Proof.} 
(i) Suppose $w_1,w_2\in F_n$ are  Reidemeister equivalent with respect to $\b(I)$.
Then,  there exists an element $w\in F_n$ such that $w_2=w^{\b(I)}w_1w^{-1}$.
Then, since 
\[w^{\g\b(f)}a_n^{\m}=(w^{\theta^{\m}})^{\b(I)\g}a_n^{\m}
=(a_n^{\m}wa_n^{-\m})^{\b(I)\g}a_n^{\m}=a_n^{\m}w^{\b(I)\g},\]
we have 
\[a_n^{\m}w_2^{\g} =a_n^{\m}w^{\b(I)\g}w_1^{\g}(w^{-1})^{\g}
=(w^{\g})^{\b(f)}(a_n^{\m}w_1^{\g})(w^{\g})^{-1},\]
which shows that $a_n^{\m}w_1^{\g}$ and $a_n^{\m}w_2^{\g}$ are Reidemeister equivalent with respect to $\b(f)$.

Conversely, suppose  $a_n^{\m}w_1^{\g}$ and $a_n^{\m}w_2^{\g}$ are Reidemeister equivalent with respect to $\b(f)$.
Then,  there exists an element $u\in F_n$ such that $a_n^{\m}w_2^{\g}=u^{\b(f)}a_n^{\m}w_1^{\g}u^{-1}$.
Let $v=u^{\g^{-1}}$.
Then, since $u^{\b(f)}a_n^{\m}=v^{\g\b(f)}a_n^{\m}=a_n^{\m}v^{\b(I)\g}$,
we have $a_n^{\m}w_2^{\g}=a_n^{\m}v^{\b(I)\g}w_1^{\g}(v^{-1})^{\g}$.
Therefore, $w_2=v^{\b(I)}w_1v^{-1}$, which shows that $w_1$ and $w_2$ are Reidemeister equivalent with respect to $\b(I)$.

(ii) This follows easily from (i).
\endproof

By (12) and Lemma 2(ii), it is enough for the proof of Theorem 1 to show that
\[\P_{\b(I)}(\tr \J(\b(I)))=\P_{\b(I)}\Bigl(\sum_{\calB\in \calP(d)}W_I(\calB)\Bigr).\]
 To prove this equality, we shall compute the matrix $ \J(\b(I))$.
 First consider the case where $I$ has length one.
For positive integers $m$, let
\begin{equation*}\G_m=\begin{cases} g_2 +\cdots +g_{m} &\text{if}\quad m\ge 2,\\
0& \text{if}\quad m=1.\end{cases}
\end{equation*}
Then, we have
\begin{Lemma} For positive integers $m$, let 
\begin{equation*}A_m=
\begin{pmatrix}
\G_m & -g_m  & -(\G_ma_2^{-1}+g_{m-1}) & 0  & \hdots & 0    \\
-a_2 &   0   &  1                     &  0  &        & \vdots\\
-a_3 &   0    &  0                     & 1   &\ddots  &\vdots\\
\vdots&\vdots&         \vdots               & 0   &\ddots & 0    \\
\vdots& \vdots     &     \vdots                   &      &\ddots & 1   \\
-a_{n-1}& 0  &  0                    &  0 &  \hdots    & 0 
\end{pmatrix}.
\end{equation*}
Then,  $\J(\b(m))=A_m^{\b(m)}$.
\end{Lemma}

{\it Proof.}
 For all $1\le i \le n-1$, we have
\begin{equation*}
a_i^{\s_1^m}=\begin{cases}
a_2^{\frac{m-1}{2}}a_2a_1^{-1}a_2^{-\frac{m-1}{2} }& \text{if}\quad  i=1, m \;\text{is odd},\\
a_2^{\frac{m}{2}}a_1a_2^{-\frac{m}{2} }& \text{if}
\quad i=1, m\; \text{is even},\\a_{i} & \text{if}\quad  2 \le i\le n-1.
\end{cases}
\end{equation*}
Also, $a_i^{\r}=a_{i+1}a_1^{-1}$.
These imply that 
\begin{equation}
a_i^{\b(m)}=\begin{cases}
(a_3a_1^{-1})^{\frac{m-1}{2}}a_3a_2^{-1}(a_3a_1^{-1})^{-\frac{m-1}{2}}
 & \text{if}\quad i=1, \;m\;\text{is odd},\\
(a_3a_1^{-1})^{\frac{m}{2}}a_2a_1^{-1}(a_3a_1^{-1})^{-\frac{m}{2} }
&  \text{if}\quad i=1, \;m \; \text{is even},\\
a_{i+1}a_1^{-1} & \text{if}\quad 2 \le i\le n-1.\end{cases}
\end{equation}
We first compute $\partial a_1^{\b(m)}/\partial a_j$ for $j=1,\ldots,n-1$.
Let $v=(a_3a_1^{-1})^{[m/2]}$,
where $[m/2]$ denotes the largest integer which does not exceed $m/2$.
Since $(a_3a_1^{-1})^r=(a_2^{\b(m)})^r=-g_{2r}^{\b(m)}$ for any positive $r$,  we have
\begin{equation}\frac{\partial v}{\partial a_j}=\begin{cases}
-\sum_{r=1}^{[m/2]}(a_3a_1^{-1})^r=\sum_{r=1}^{[m/2]}g_{2r}^{\b(m)} & \text{if}\quad j=1,\\
\sum_{r=0}^{[m/2]-1}(a_3a_1^{-1})^r=-\displaystyle{\frac{\partial v}{\partial a_1}}
(a_2^{-1})^{\b(m)}& \text{if} \quad j=3, \\
0 & \text{otherwise}.
\end{cases}
\end{equation}
It follows from the definition of $g_{2r}$ that 
\begin{equation}(1-a_1)\sum_{r=1}^{[m/2]}g_{2r}=\begin{cases}
\G_{m} & \text{if} \quad m \;\text{is odd,}\\
\G_{m+1} & \text{if} \quad m \;\text{is even.}
\end{cases}
\end{equation}

Suppose $m$ is odd. Let $w=a_3a_2^{-1}$. Then $vwv^{-1}=a_1^{\b(m)}$ by (13).
We shall compute the right-hand side of the equality (7).
We have by (14), (15)
\begin{equation*}
(1-vwv^{-1})\frac{\partial v}{\partial a_j}=\begin{cases}
\G_{m}^{\b(m)} & \text{if}\quad j=1,\\
-(\G_{m}a_2^{-1})^{\b(m)} & \text{if}\quad j=3,\\
0 & \text{otherwise},
\end{cases}
\end{equation*}
and, since $v=(a_2^{[m/2]})^{\b(m)}$, we have 
\begin{equation*}
v\frac{\partial w}{\partial a_j}=\begin{cases}
-vw =-a_1^{\b(m)}v=-g_m^{\b(m)} & \text{if}\quad j=2,\\
v=-g_{m-1}^{\b(m)} & \text{if}\quad j=3,\\
0 & \text{otherwise}.
\end{cases}
\end{equation*}

Suppose $m$ is even. Let $w'=a_2a_1^{-1}$. Then $vw'v^{-1}=a_1^{\b(m)}$ by (13). 
By (14), (15), we have 
\begin{equation*}
(1-vw'v^{-1})\frac{\partial v}{\partial a_j}=\begin{cases}
\G_{m+1}^{\b(m)} & \text{if}\quad j=1,\\
-(\G_{m+1}a_2^{-1})^{\b(m)}=-(\G_{m}a_2^{-1}+g_{m-1})^{\b(m)} & \text{if}\quad j=3,\\
0 & \text{otherwise}.
\end{cases}
\end{equation*}

\begin{equation*}
v\frac{\partial w'}{\partial a_j}=\begin{cases}
-vw' =-a_1^{\b(m)}v=-g_{m+1}^{\b(m)}& \text{if}\quad j=1,\\
v=-g_m^{\b(m)} & \text{if}\quad j=2,\\
0 & \text{otherwise}.
\end{cases}
\end{equation*}
These computations and the equality (7)  imply that 
$\partial a_1^{\b(m)}/\partial a_j$ is equal to the $(1,j)$ entry of the matrix $A_m^{\b(m)}$
in either case of $m$ odd or even.

The $i$-th row for  $i\ge 2$ is obtained from the following: 
 
\[\frac{\partial a_i^{\b(m)}}{\partial a_j}=\frac{\partial a_{i+1}a_1^{-1}}{\partial a_j}=
\begin{cases}
-a_i^{\b(m)}& \text{if}\quad j=1,\\
1& \text{if}\quad j=i+1,\\
0& \text{otherwise}.\end{cases}\]
This completes the proof.
\endproof

Now consider the general case of $I$ having an arbitrary length.
Fix a sequence $I=(i_1,\ldots,i_d)$ of positive integers.
We shall give  formulas for entries of the matrix $\J(\b(I))$ in Lemma 5 below. To state these formulas, we need four families 
of elements of $\Z F_n$. The first one is $\a_q^{l}$ defined for  integers $q,l$.
For integers $1\le l\le d$, denote  $\b_l(I)$ simply by $\b_l$, and let $\b_{d+1}=e$.
Then,  $\a_q^{l}\in \Z F_n$ is  defined for integers $q,l$ by
\[\a_q^{l}=\begin{cases} a_q^{\b_{l}} & \text{if} \quad 1\le q\le n-1,1\le l\le d,\\
0 & \text{otherwise}.\end{cases}\] 

The second family is  $W_k^l\in \Z F_n$ defined for positive integers $k,l$. To define these elements, we need to generalize the notion of a partition given in 
Definition 2 to the non-cyclic setting.
 
\begin{Definition} Suppose    $k,l$ are positive integers with $k\le l$. 
\begin{enumerate}
\item  For integers $p,q$ with $k\le p\le q\le l$, define a sequence $[p,q]$ of positive integers 
by $[p,q]=(p, \ldots,q)$. This sequence is called  a {\it block} in $\{k,\ldots,l\}$, and 
the number of integers contained in it is called its {\it length}.
For a block $B$, let $\underline{B}$ denote its underlying set.
\item A set $\{B_1,\ldots,B_s\}$ of blocks in $\{k,\ldots,l\}$ is a {\it partition} of $\{k,\ldots,l\}$
if $\underline{B_1},\ldots,\underline{B_s}$ are mutually disjoint, 
$\underline{B_1}\cup\cdots \cup \underline{B_s}=\{k,\ldots,l\}$, and  $B_1,\ldots,B_s$ have length less than or equal to $n-1$.
\item Let $\calP(k,l)$ denote the set of partitions of $\{k,\ldots,l\}$.
\end{enumerate}
For a subset  $\calA$ of $\calP(d)$ or of $\calP(k,l)$, where $1\le k\le l\le d$,
let $W_I(\calA)=\sum_{\calB\in \calA}W_I(\calB)$.
Then,  $W_k^l\in \Z F_n$ is defined for positive integers $k,l$ by 
 \[W_k^l=\begin{cases} W_I(\calP(k,l)) & \text{if}\quad k\le l\le d,\\
1 & \text{if}\quad k=l+1\;\text{and}\; l\le d\\
0 & \text{otherwise}.
\end{cases}\]
\end{Definition}
We prepare the next lemma, which will be used to prove Lemma 5.
For $1\le l \le  d$, $1\le \l\le n-1$,
 let $\calP_{1,\l}(l)$ be the set of partitions of $\{1,\ldots,l\}$ 
such that the block with initial element $1$ has length $\l$, and
 let $\calP_{l,\l}(d)$ be the set of partitions of $\Z_d$ 
which contain a block  with initial element $l$ and length $\l$.

\begin{Lemma} 
\begin{enumerate}
\item $g_m^{\b(m)}=-g_{m-1}^{\b(m)}a_3a_2^{-1}$ for any  positive integer  $m$.
\item For  $1\le l\le d$, we have
\[\sum_{u\ge 0} \left( g_{i_1}^{\b_{1}}\a_{2+u}^{2+u}+g_{i_1-1}^{\b_1}\a_{3+u}^{2+u}\right)W_{3+u}^{l}=
 W_I(\calP_{1,n-1}(l)).\]
\item For positive integers $l$ with  $d+3-n \le l\le d$, the elements
\[\sum_{u\ge 0} \left(\a_{d+2-l+u}^{1+u}W_{2+u}^{l-1}g_{i_l}^{\b_l}+
\a_{d+3-l+u}^{1+u}W_{2+u}^{l-1}g_{i_{l}-1}^{\b_{l}}\right)\]
and   $W_I(\calP_{l,n-1}(d))$ have the same $\P_{\b(I)}$-image.
\end{enumerate}
\end{Lemma}

{\it Proof}. (i) Consider the case of  $m$ odd.
Since $a_3a_1^{-1}=a_2^{\b(m)}$, we have  
\[(a_3a_1^{-1})^{(m-1)/2}= (a_2^{\b(m)})^{(m-1)/2}=(a_2^{(m-1)/2})^{\b(m)}=c_{m-1}^{\b(m)}=-g_{m-1}^{\b(m)}.\]
Therefore, we have by (13), $a_1^{\b(m)}=-g_{m-1}^{\b(m)}a_3a_2^{-1}(c_{m-1}^{\b(m)})^{-1}$ and hence 
\[g_m^{\b(m)}=(a_1a_2^{(m-1)/2})^{\b(m)}=a_1^{\b(m)}c_{m-1}^{\b(m)}=-g_{m-1}^{\b(m)}a_3a_2^{-1}.\]
In the case of  $m$ even, we have by (13)
 the desired equality from the following: 
\begin{align*}
g_{m-1}^{\b(m)} &=(a_1a_2^{m/2}a_2^{-1})^{\b(m)}\\
&=\left[ c_m^{\b(m)}a_2a_1^{-1}(c_m^{\b(m)})^{-1}\right] c_m^{\b(m)}(a_2^{-1})^{\b(m)}\\
&=c_m^{\b(m)}a_2a_1^{-1}(a_3a_1^{-1})^{-1}\\
&=c_m^{\b(m)}a_2a_3^{-1}=-g_m^{\b(m)}a_2a_3^{-1}.
\end{align*}
(ii) Let 
$\Sigma_1$  be the left-hand side of the equality (ii). 
Let 
\[V(u)=g_{i_1}^{\b_1}\a_{2+u}^{2+u}W_{3+u}^l,\quad V'(u)=g_{i_1-1}^{\b_1}\a_{3+u}^{2+u}W_{3+u}^l.\]
Then, we have \[\Sigma_1=\sum_{u\ge 0}(V(u)+V'(u)).\]
There are three cases:
\begin{enumerate}[(a)]
\item$2+u> l$ or $2+u>n-1$, 
\item$2+u\le l$ and $2+u<n-1$,
\item $2+u\le l$ and $2+u=n-1$.
\end{enumerate}

Consider Case (a). If $2+u>l$, then $3+u>l+1$ and so $W_{3+u}^l=0$.
Also, if $2+u>n-1$, then $\a_{2+u}^{2+u}=\a_{3+u}^{2+u}=0$. Therefore,  we have   $V(u)=V'(u)=0$.
Consider Case (b). 
For any $i\ge 2$, $1\le l\le d$, and any positive integer $u$ with $l+u\le d+1$ and $i+1+u\le n-1$, 
we have by the equality $a_i^{\b(m)}=a_{i+1}a_1^{-1}$ for any $m$ that
\begin{equation}(a_{i+1}a_i^{-1})^{\b_l}
=(a_{i+1}a_i^{-1})^{\b(i_{l},\ldots,i_{l+u-1})\b_{l+u}}
=(a_{i+1+u}a_{i+u}^{-1})^{\b_{l+u}}.\end{equation}
We have by (i) of this lemma and (16) that  
\begin{align*} 
g_{i_1}^{\b_1}&=(g_{i_1}^{\b(i_1)})^{\b_2}
=-(g_{i_1-1}^{\b(i_1)}a_3a_2^{-1})^{\b_2}\\
&=-g_{i_1-1}^{\b_1}a_{3+u}^{\b_{2+u}}(a_{2+u}^{\b_{2+u}})^{-1}.
\end{align*}
This implies that  $V(u)+V'(u)=0$. Consider Case (c).
Since $3+u=n$, we have $\a_{3+u}^{2+u}=0$ and hence $V'(u)=0$.
Also, since  $n-1=2+u\le l\le d$ and hence  $\a_{2+u}^{2+u}=a_{n-1}^{\b_{n-1}}$, we have 
 \[V(u)=g_{i_1}^{\b_1}a_{n-1}^{\b_{n-1}}W_n^{l}
=W_I([1,n-1])W_n^l=W_I(\calP_{1,n-1}(l)).\]
If $l\ge n-1$,  putting these computations together, we have
 $\Sigma_1=W_I(\calP_{1,n-1}(l))$, and
so (ii) holds. Suppose  $l<n-1$.
Then, $\Sigma_1=0$ since Case (c) does not occur, and we have  $W_I(\calP_{1,n-1}(l))=0$
 since $\calP_{1,n-1}(l)$ is empty. 
Therefore, the  equality (ii) is proved.

(iii) Let $l$ be a positive integer with $d+3-n\le l\le d$. Note that $d+2-l\le n-1$.
Let
\[V_{l}(u)=\a_{d+2-l+u}^{1+u}W_{2+u}^{l-1}g_{i_l}^{\b_l},\quad
 V'_l(u)=\a_{d+3-l+u}^{1+u}W_{2+u}^{l-1}g_{i_{l-1}}^{\b_{l}},\]
and let \[\Sigma_l=\sum_{u\ge 0}(V_l(u)+V'_l(u)).\]
There are three cases:

\begin{enumerate}[(a)]
\item  $2+u> l$ or $d+2-l+u>n-1$, 
\item $2+u\le l$ and $d+2-l+u<n-1$,
\item $2+u\le l$ and $d+2-l+u=n-1$.
\end{enumerate}

Consider Case (a). If $ 2+u>l$, then $W_{2+u}^{l-1}=0$. Also, if
$d+2-l+u>n-1$, then $\a_{d+2-l+u}^{1+u}=\a_{d+3-l+u}^{1+u}=0$.
Therefore, we have  $V_l(u)=V'_l(u)=0$. In Case (b), we have by (i) of this lemma and (16) that
\begin{align*}
g_{i_l}^{\b_l}&=(g_{i_l}^{\b(i_l)})^{\b_{l+1}}
=-g_{i_l-1}^{\b_l}(a_3a_2^{-1})^{\b_{l+1}\b_1\b_1^{-1}}\\
&=-g_{i_l-1}^{\b_l}(a_{3+d-l}a_{2+d-l}^{-1})^{\b_1\b_1^{-1}}
=-g_{i_l-1}^{\b_l}(a_{d+3-l+u}^{\b_{1+u}}(a_{d+2-l+u}^{\b_{1+u}})^{-1})^{\b_1^{-1}}.
\end{align*}
Therefore, noting that $\b_1=\b(I)$, we have by (5),
\begin{align*}
\P_{\b(I)}(V_l(u))&=\P_{\b(I)}(W_{2+u}^{l-1}g_{i_l}^{\b_l}(a_{d+2-l+u}^{\b_{1+u}})^{\b_1^{-1}})\\
&=-\P_{\b(I)}(W_{2+u}^{l-1}g_{i_l-1}^{\b_l}(a_{d+3-l+u}^{\b_{1+u}})^{\b_1^{-1}})=-\P_{\b(I)}(V'_l(u)).
\end{align*}
Therefore, $\P_{\b(I)}(V_l(u)+V'_l(u))=0$.
In Case (c), clearly $V'_l(u)=0$.
Since $\a_{d+2-l+u}^{1+u}=a_{n-1}^{\b_{1+u}}$ and 
 the length of the block $[l,1+u]$ is $(d-l+1)+1+u=n-1$,
we have 
 \begin{align*}
\P_{\b(I)}(V_l(u))&=\P_{\b(I)}(a_{n-1}^{\b_{1+u}}W_{2+u}^{l-1}g_{i_l}^{\b_l})
=\P_{\b(I)}(W_{2+u}^{l-1}g_{i_l}^{\b_l}a_{n-1}^{\b_{1+u}\b_1^{-1}})\\
&=\P_{\b(I)}(W_{2+u}^{l-1}W_I([l,1+u]))\\
&=\P_{\b(I)}(W_I(\calP_{l,n-1}(d))).
\end{align*}
If $d<n-1$, Case (c) does not occur since $d+2-l+u=n-1$ implies $2+u=(n-1)-d+l>l$.
Therefore, (iii) is proved by summing up these computations.
\endproof

The last two families of elements of $\Z F_n$ necessary to state Lemma 5 are $S_l,G_l$ 
defined for integers  $l$ as follows:
\begin{equation*}
S_l=\begin{cases}
     (\G_{i_l}a_2^{-1}+g_{i_l-1})^{\b_l} & \text{if}\quad 1\le l\le d,\\
      0 & \text{otherwise,}
   \end{cases}
\end{equation*}
\begin{equation*}
G_l=\begin{cases}
      g_{i_l}^{\b_l} & \text{if}\quad 1\le l\le d,\\
      -1 & \text{if} \quad l=d+1,\\
      0 & \text{otherwise.}
   \end{cases}
\end{equation*}
Let $r_{i,j}(I)$ be the $(i,j)$-entry of the matrix $\J(\b(I))$.
\begin{Lemma}
\begin{equation*}
r_{i,j}(I)=
\begin{cases}
-W_1^{d+2-j}S_{d+3-j}-W_1^{d+1-j}G_{d+2-j}& \text{if} \quad i=1,\\
\sum_{u\ge 0}\a_{i+u}^{1+u}\left (W_{2+u}^{d+2-j}S_{d+3-j}+W_{2+u}^{d+1-j}G_{d+2-j}\right) 
+\d_{i,j-d} & \text{if} \quad i\ge 2.
\end{cases}
\end{equation*}
\end{Lemma}

{\it Proof}. We prove this lemma by the induction on $d$.
The case of $d=1$ follows easily from Lemma 3.
Assume that the lemma holds for $d-1$, and we shall prove it for $d$.
 Let $I=(i_1,\ldots,i_{d})$ be a sequence of positive integers.
Let  $I'=(i_2,\ldots,i_d)$.
Then, $r_{i,j}(I')$ is obtained from the right-hand side of this lemma by replacing 
 $\a_{i+u}^{1+u}$ and $\d_{i,j-d}$ with $\a_{i+u}^{2+u}$ and $\d_{i,j-(d-1)}$ respectively,
and by adding one to the subscript of each of $W$'s.
Note that by (10) and Lemma 3
\begin{equation}\J(\b(I))=\J(\b(i_1))^{\b(I')}\J(\b(I'))=A_{i_1}^{\b_1}\J(\b(I')).\end{equation}

Consider the case of $i=1$. 
Let
\[M(l)=\G_{i_1}^{\b_1}W_2^{l}+(\G_{i_1}a_2^{-1})^{\b_1} \sum_{u\ge 0}\a_{3+u}^{2+u}W_{3+u}^{l}
+W_I(\calP_{1,n-1}(l))\]
for $l\ge 1$ and $M(l)=0$ for $l\le 0$.
Then, we have by (17), Lemma 3, and Lemma 4(ii)  that $r_{1,j}(I)$ is equal to 
\begin{align*}
&\G_{i_1}^{\b_1}r_{1,j}(I')-g_{i_1}^{\b_1}r_{2,j}(I')-(\G_{i_1}a_2^{-1}
+g_{i_1-1})^{\b_1}r_{3,j}(I')\\
=&-\G_{i_1}^{\b_1}\left(W_2^{d+2-j}S_{d+3-j}+W_2^{d+1-j}G_{d+2-j}\right)\\
&   -g_{i_1}^{\b_1} \left[\sum_{u\ge 0}\left( \a_{2+u}^{2+u}W_{3+u}^{d+2-j}S_{d+3-j}
+\a_{2+u}^{2+u}W_{3+u}^{d+1-j}G_{d+2-j}\right) +\d_{2,j-(d-1)}\right]\\
&  -(\G_{i_1}a_2^{-1}+g_{i_1-1})^{\b_1}
\left[\sum_{u\ge 0}\left( \a_{3+u}^{2+u}W_{3+u}^{d+2-j}S_{d+3-j}
+\a_{3+u}^{2+u}W_{3+u}^{d+1-j}G_{d+2-j}\right) +\d_{3,j-(d-1)}\right]\\
=&-(M(d+2-j)S_{d+3-j}+\d_{j,d+2}S_1)-(M(d+1-j)G_{d+2-j}+\d_{j,d+1}G_1).
\end{align*}
Therefore, since $\d_{j,d+2}S_1=\d_{j,d+2}S_{d+3-j}$ and $\d_{j,d+1}G_1=\d_{j,d+1}G_{d+2-j}$,
we have \[r_{1,j}(I)=-(M(d+2-j)+\d_{j,d+2})S_{d+3-j}-(M(d+1-j)+\d_{j,d+1})G_{d+2-j}.\]
Since $M(l)=\sum_{\l=1}^{n-1}W_I(\calP_{1,\l}(l))=W_1^{l}$ if $l\ge 1$ and $W_1^0=1$,
  this is equal to $-W_1^{d+2-j}S_{d+3-j}-W_1^{d+1-j}G_{d+2-j}$, which is 
the right-hand side of the equality of the lemma in the case of $i=1$. 

Consider the case of $i\ge 2$. We have
\begin{align*}
r_{i,j}(I)
&=-a_i^{\b_1}r_{1,j}(I')+r_{i+1,j}(I')\\
&=a_i^{\b_1}({W}_2^{d+2-j}S_{d+3-j}+{W}_2^{d+1-j}G_{d+2-j})\\
& \hspace{5mm}  +\sum_{u\ge 0}\a_{i+1+u}^{2+u}({W}_{3+u}^{d+2-j}S_{d+3-j}
+{W}_{3+u}^{d+1-j}G_{d+2-j})+\d_{i+1,j-(d-1)}.
\end{align*}
It is easy to show that this is equal to the right-hand side of the equlatity  in  the lemma 
in the case of $i\ge 2$. 
\endproof

We shall complete the proof of Theorem 1.
Note that $S_{d+1}=S_{d+2}=0$. Let
\[L_1=\sum_{j=3}^{\n_1}\sum_{u\ge 0}\a_{j+u}^{1+u}W_{2+u}^{d+2-j}(\G_{i_{d+3-j}}a_2^{-1})^{\b_{d+3-j}},\]
\[L_2=\sum_{j=3}^{\n_2}\sum_{u\ge 0}\a_{j+u}^{1+u}W_{2+u}^{d+2-j}g_{i_{d+3-j}-1}^{\b_{d+3-j}},\]
 \[L_3=\sum_{j=2}^{\n_3}\sum_{u\ge 0}\a_{j+u}^{1+u}W_{2+u}^{d+1-j}g_{i_{d+2-j}}^{\b_{d+2-j}},\]
where $\n_1=\n_2=\min\{n-1,d+1\},\n_3=\min\{n-1, d\}$.
Then, by Lemma 5,
\[\tr(\J(\b(I)))= W_1^{d}+L_1+L_2+L_3.\]
Let $d_1=d+3-\n_1=\max\{d+4-n,2\}, d_2=d+2-\n_3=\max\{d+3-n,2\}$.
In $L_2$, we can change $\n_2$ to $\bar{\n}_2=\min\{n,d+1\}$, since $\a_{n+u}^{1+u}=0$.
Then, $d+3-\bar{\n}_2=d_2$.
Putting $l =d+3-j$ in $L_1$ and $L_2$, and putting $l=d+2-j$ in $L_3$, we have
\[L_1=\sum_{l=d_1}^d
\sum_{u\ge 0}\a_{d+3-l+u}^{1+u}W_{2+u}^{l-1}(\G_{i_{l}}a_2^{-1})^{\b_{l}},\]
\[L_2+L_3=\sum_{l=d_2}^{d}\sum_{u\ge 0}\left(\a_{d+3-l+u}^{1+u}W_{2+u}^{l-1}g_{i_{l}-1}^{\b_l}
+\a_{d+2-l+u}^{1+u}W_{2+u}^{l-1}g_{i_{l}}^{\b_l}\right).\]
Let \[Q_1=\bigcup_{l=d_1}^d\bigcup_{\l=d+2-l}^{n-2}\calP_{l,\l}(d),\quad
Q_2=\bigcup_{l=d_2}^d\calP_{l,n-1}(d).\]
Since $\P_{\b(I)}(\a_{d+3-l+u}^{1+u}W_{2+u}^{l-1}(\G_{i_{l}}a_2^{-1})^{\b_{l}})$
is equal to \[\P_{\b(I)}(W_{2+u}^{l-1}W_I([l,1+u]))=\P_{\b(I)}(W_I(\calP_{l,d+2-l+u}(d)))\]
if $d+2-l+u\le n-2$, and equal to zero otherwise since $\a_{m}^{1+u}=0$ for $m\ge n$,
we have 
\[\P_{\b(I)}(L_1)=\sum_{l=d_1}^{d}\sum_{\l=d+2-l}^{n-2}\P_{\b(I)}(W_I(\calP_{l,\l}(d)))
=\P_{\b(I)}(W_I(Q_1)).\]
Lemma 4(iii) implies that \[\P_{\b(I)}(L_2+L_3)
=\sum_{l=d_2}^d \P_{\b(I)}(W_I(\calP_{l,n-1}(d)))=\P_{\b(I)}(W_I(Q_2)).\]
Since  $W_1^{d}=W_I(\calP(1,d))$ and  $\calP(d)=\calP(1,d)\cup Q_1\cup Q_2$,
these equalities prove that  $\P_{\b(I)}(\tr(\J(\b(I)))$ 
is equal to $\P_{\b(I)}(W_I(\calP(d)))$.
Thus the proof of the theorem is completed by (12) and Lemma 2(ii).

\section{Reduction of the formula}

This section makes  preparations for the proof of Proposition 2.
We shall show that, under the assumption of Proposition 2, 
the element $\sum_{\calB\in \calP(d)}W_I(\calB)$ of $\Z F_n$ in the formula of 
Theorem 1 can be  reduced so that 
Reidemeister equivalent elements of $F_n$ have the same coefficient.
Hence, no cancellation occurs when the reduced one is projected on $\Z\calR(\b(f))$,
which enables us to apply Lemma 1 to the problem.

Consider first the case where  $n\ge 4$ and $i_1,\ldots,i_d\ge 2$.
Let $\calP'(d)$ be the set of partitions $\calB=\{B_1,\ldots, B_s\}$ of $\Z_d$ such that
$(|B_j|,|B_{j+1}|)\ne (1,n-2)$ for any $1\le j\le s$, where $B_{s+1}=B_1$.
For partitions $\calB\in \calP'(d)$, we shall define elements $W_I'(\calB)$ of $\Z F_n$ which satisfy that  the sums $\sum_{\calB\in\calP(d)}W_I(\calB)$ and $\sum_{\calB\in\calP'(d)}W_I'(\calB)$ have the same image under the projection
$\P_{\b(I)}$.
Suppose  $B=[p,q]$ is a block. If  $|B|<n-1$, let $S_B(I)$ 
denote the set of integers $J$ with $0\le J\le i_p-2$, and
let $\l_B(J)=c_{J}^{\a(B)}a_{|B|+1}^{\o(B)}\in F_n$ for any $J\in S_B(I)$. 
If $|B|=n-1$,  let $S_B(I)$ denote the set of $(j,j')\in \Z^2$ such that 
 $2\le j\le i_p, \;\;0\le j'\le i_{p'}-2, \;\;(j,j')\ne (i_p,0)$, 
and let $\l_B(J)= c_j^{\a(B)}c_{j'}^{\a'(B)}a_{n-1}^{\o(B)}\in F_n$ for $J=(j,j')\in S_B(I)$,
where $\a'(B)\in B_n$ is defined by 
\begin{equation}
\a'(B)=\begin{cases}
     \b_{p+1} & \quad \text{if}\; p\le d-1, \notag \\
   e &  \quad\text{if}\;p=d. \notag\end{cases}
\end{equation} 
 For a partition $\calB=\{B_1,\ldots,B_s\}$, let 
\[S_{\calB}(I)=S_{B_1}(I)\times \cdots \times S_{B_s}(I).\]
For an element $\calJ=(J_1,\ldots,J_s)$ of $S_{\calB}(I)$, define $\l_{\calB}(\calJ)\in F_n$ by  
$\l_{\calB}(\calJ)=\l_{B_1}(J_1)\cdots \l_{B_s}(J_s)$. 
For a block $B=[p,q]$, define $W'_I(B)\in \Z F_n$ by
\begin{equation}
W'_I(B)=\begin{cases}
W_I(B)=\displaystyle{\sum_{J\in S_B(I)}g_J^{\a(B)}a_{|B|+1}^{\o(B)}} & \text{if}\quad |B|< n-1,\notag \\
\displaystyle{\sum_{(j,j')\in S_B(I)} }g_j^{\a(B)}g_{j'}^{\a'(B)}a_{n-1}^{\o(B)} & \text{if} \quad |B|= n-1.
\end{cases}
\end{equation} 
Then, 
$W_I'(\calB)\in \Z F_n$ is defined for $\calB\in \calP'(d)$ by  $W_I'(\calB)= W_I'(B_1)\cdots W_I'(B_s)$,
where $\calB=\{B_1,\ldots,B_s\}$ with $1\le p_1<\cdots <p_s\le d$.

For $w\in F_n$, define an integer  $e(w)$ as the exponent sum of $w$ 
with respect to the standard generators $\xi_1,\ldots,\xi_n$. Note that $e(w)$
can be defined also by  $\calT(w)=t^{e(w)}$, where $\calT$ is the ring 
homomorphism introduced in Example 4.

\begin{Lemma} Let  $n\ge 4$.
Assume $\b(f)=\g^{-1}\theta^{\m}\b(I)\g$, where $\m\in\Z, \g\in B_n$,
 $I=(i_1,\ldots,i_d)$ with  $i_1,\ldots,i_d \ge 2$. Then,  we have 
 \begin{enumerate}
 \item \[\calL(f)= -\P_{\b(f)}\Bigl(a_n^{\m}\sum_{\calB\in \calP'(d)}W'_I(\calB)^{\g}\Bigr).\]
\item For $\calB\in \calP'(d)$, we have  \[W_I'(\calB)=
\sum_{\calJ\in S_{\calB}(I)}(-1)^{d+e(\l_{\calB}(\calJ))}\l_{\calB}(\calJ).\]
\item For any $\calB\in \calP'(d)$ and 
any $\calJ\in S_{\calB}(I)$, the coefficient of
$\P_{\b(f)}(a_n^{\m}\l_{\calB}(\calJ)^{\g})$ in $\calL(f)$ is non-zero.
\end{enumerate}
\end{Lemma}

{\it Proof}. (i) For a partition $\calB$, let 
$K(\calB)$ be the set of integers $k\in \Z_d$  such that either 
$[k,k+n-2]$ is a block in $\calB$, or both  $(k)$ and $[k+1,k+n-2]$ are blocks in $\calB$,
where integers are taken modulo $d$. 
Let $\calK(d)$ be the set of subsets $K$ of $\Z_d$ such that, if 
 $K$ is written as $K=\{k_1,\ldots,k_t\}$, where $1\le k_1<\cdots <k_t\le d$, then
$k_{r+1}-k_r\ge n-1$ for any $1\le r \le t$, where
we put $k_{t+1}=k_1+d$. We assume that the empty set is contained in $\calK(d)$.
Note that a subset $K$ of $\Z_d$ is contained in $\calK(d)$
if and only if there is a partition $\calB$ with $K(\calB)=K$.
For  $K\in\calK(d)$, 
let $\calP_K$ be the set of partitions $\calB$ with $K(\calB)=K$. 

Assume $d\ge n-1$. Let $K=\{k_1,\ldots,k_t\}\in \calK(d)$, where $k_1<\cdots <k_t$.
For $1\le r \le t$, let $B(r)=[k_r,k_r+n-2]$ and 
let $X_r=W_I(B(r))\in\Z F_n$.
Also, define $Y_r\in \Z F_n$ by
\[Y_r=\begin{cases}W_I((d))W_I([1,n-2])^{\b_1^{-1}} & \text{if} \;\; r=t, k_t=d,\\
W_I((k_r))W_I([k_{r}+1,k_{r}+n-2])& \text{otherwise}.\end{cases}\] 
For $1\le k,l\le d$ with $k\le l+1$, define $Z(k,l)\in \Z F_n$ by
 $Z(k,l)=W_I(\calP(k,l))$ if $k\le l$ and $Z(k,l)=e$ if $k=l+1$.
For $r$ with $1\le r<t$,
let  $Z_r=Z(k_r+n-1, k_{r+1}-1)$. 
If  $k_t+n-2<d$, let $Z_0=Z(1,k_1-1)$ and $Z_t=Z(k_t+n-1,d)$.
If   $k_t+n-2\ge d$, let $Z_0=Z(k_t+n-1-d,k_1-1)$ and $Z_t=e$.
Let 
\[\L_K'=\prod_{r=1}^{t-1} (X_r+Y_r)Z_r,\quad \L_K=Z_0 \L_K'(X_t+Y_t)Z_t.\]
In the case of $k_t<d$, it is easy to see that $\sum_{\calB\in \calP_K}W_I(\calB)=\L_K$.
In the case of $k_t=d$, we have  $\sum_{\calB\in \calP_K}W_I(\calB)=Z_0\L_K'X_t+W_I([1,n-2])Z_0\L_K'W_I((d))$.  By (5), this has the same image as $\L_K$ under 
the projection $\P_{\b(I)}$.
Therefore, in either case, we have 
\begin{equation} \textstyle{\P_{\b(I)}\bigl(\sum_{\calB\in \calP_K}W_I(\calB)\bigr)}=\P_{\b(I)}(\L_K).
\end{equation}
Note that  $W_I(B)=W_I'(B)$ for any block  $B$ in $\calB\in \calP_K$ with $\underline{B}$ 
 disjoint from  $\underline{B(1)}\cup\cdots \cup \underline{B(t)}$.
 Also, letting $\a_r=\a(B(r)),\a'_r=\a'(B(r))$ and $\o_r=\o(B(r))$, we have   
\begin{align*}X_r+Y_r&=g_{i_{k_r}}^{\a_r}a_{n-1}^{\o_r}
+\G_{i_{k_r}}^{\a_r}(\G_{i_{k_{r}+1}}a_2^{-1})^{\a_r'}a_{n-1}^{\o_r} \notag\\
&=\sum_{(j,j')\in S_{B(r)}(I)} g_{j}^{\a_r}g_{j'}^{\a_r'}a_{n-1}^{\o_r}\notag \\
&=W_I'(B(r)).
\end{align*}
Hence, we have the equality \[\L_K=\sum_{\calB\in\calP_K\cap\calP'(d)}W_I'(\calB).\]
This and (18) imply that $\sum_{\calB\in\calP_K}W_I(\calB)$ 
and $\sum_{\calB\in\calP_K\cap\calP'(d)}W_I'(\calB)$ have the same image 
under $\P_{\b(I)}$.
Furthermore, since the disjoint unions $\cup_{K\in\calK(d)}\calP_K$ and $\cup_{K\in\calK(d)}(\calP_K\cap\calP'(d))$ coincide with  $\calP(d)$ and $\calP'(d)$ respectively,
we have \[\sum_{\calB\in\calP(d)}W_I(\calB)=\sum_{K\in \calK(d)}\sum_{\calB\in \calP_K}W_I(\calB),\quad \sum_{K\in \calK(d)}\sum_{\calB\in\calP_K\cap\calP'(d)}W_I'(\calB)=\sum_{\calB\in\calP'(d)}W_I'(\calB).\]
Therefore, $\sum_{\calB\in\calP(d)}W_I(\calB)$ and 
$\sum_{\calB\in\calP'(d)}W_I'(\calB)$ have the same image under $\P_{\b(I)}$,
and  (i) follows from Theorem 1. 

In the case of $d<n-1$, it is trivial that $\calP'(d)=\calP(d)$ 
and $W_I'(\calB)=W_I(\calB)$  for any partition $\calB$.
Hence, the formula (i) is identical with that in Theorem 1. 

(ii) Let $B$ be a block. Then, by the definition of $\l_{B}(J)$, we see that
$W_I'(B)$ is written in  the form $W_I'(B)=\sum_{J\in S_B(I)}\e(J)\l_B(J)$, where $\e(J)$
are integers. We have 
\begin{equation} \e(J)=(-1)^{|B|+e(\l_B(J))}.\end{equation}
In fact, if $|B|<n-1$,  we have $\e(J)=(-1)^{J+1}$ and this is
 equal to $(-1)^{|B|+e(\l_B(J))}$ since $e(\l_B(J))=J+|B|+1$.  
Also, if $|B|=n-1$, $\e(J)=(-1)^{j+j'}$ and this is equal to $(-1)^{|B|+e(\l_B(J))}$ 
since $e(\l_B(J))=j+j'+n-1=j+j'+|B|$.
Let $\calB=\{B_1,\ldots,B_s\}\in \calP'(d)$ and $\calJ=\{J_1,\ldots,J_s\}\in S_{\calB}(I)$.
Then, since $|B_1|+\cdots +|B_s|=d$ and 
 $e(\l_{B_1}(J_1))+\cdots +e(\l_{B_s}(J_s))=e(\l_{\calB}(\calJ))$, 
the coefficient of $\l_{\calB}(\calJ)$ in $W_I'(\calB)$ is equal to   $\e(J_1)\cdots \e(J_s)$, which is equal to $(-1)^{d+e(\l_{\calB}(\calJ))}$ by (19).
Thus,  (ii) is proved.

(iii) Let $\G(I)$ be the set of pairs $(\calB,\calJ)$ with $\calB\in \calP'(d)$, $\calJ\in S_{\calB}(I)$.
We say two elements $(\calB,\calJ)$, $(\calB',\calJ')\in \G(I)$ are equivalent
if  $\l_{\calB}(\calJ)$ is Reidemeister equivalent to $\l_{\calB'}(\calJ')$ with respect to $\b(I)$.
This defines an equivalence relation on $\G(I)$. 
Denote by $[(\calB,\calJ)]$
 the equivalence class represented by $(\calB,\calJ)$.
Let $n(\calB,\calJ)$ be the coefficient of
$\P_{\b(f)}(a_n^{\m}\l_{\calB}(\calJ)^{\g})$ in $-\calL(f)$.
Then, by (i), (ii) of this lemma and  Lemma 2(i), 
$n(\calB,\calJ)$ is equal to 
the sum of the coefficient of  $\l_{\calB'}(\calJ')$ in $W_I'(\calB')$
taken over the elements $(\calB',\calJ')$  of  $[(\calB,\calJ)]$.
Since this coefficient is equal to $(-1)^{d+e(\l_{\calB'}(\calJ'))}$ by (ii),
we have \begin{equation}
n(\calB,\calJ)=\sum_{(\calB',\calJ')\in [(\calB,\calJ)]}(-1)^{d+e(\l_{\calB'}(\calJ'))}.
\end{equation}
For any   $(\calB',\calJ') \in [(\calB,\calJ)]$, we have $e(\l_{\calB'}(\calJ'))=e(\l_{\calB}(\calJ))$,
 since $\l_{\calB'}(\calJ')$ is Reidemeister equivalent to $\l_{\calB}(\calJ)$
and the exponent sum of an element of $F_n$ is preserved under the action of $B_n$ on $F_n$.
 Therefore, (20) implies that $n(\calB,\calJ)$ is equal to
$(-1)^{d+e(\l_{\calB}(\calJ))}\sharp  [(\calB,\calJ)]$, which is clearly non-zero.
\endproof

 Consider next the case where $n=3$ and $i_1,\ldots,i_d\ge 3$. 
Let $\Z^d(I)$ denote the set of $J=(j_1,\ldots,j_d)\in \Z^d$ which satisfy
$2\le j_l\le i_l$ for any $1\le l \le d$, and let $S(I)$ be the set of $J=(j_1,\ldots,j_d)\in \Z^d(I)$ 
with $(j_l,j_{l+1})\ne (i_l,2)$ for any $1\le l \le d$, where $j_{d+1}=j_1$.
For $J=(j_1,\ldots,j_d)\in \Z^d$, 
 let $|J|=j_1+\cdots +j_d$, $c(J)=c_{j_1}^{\b_1}\cdots c_{j_d}^{\b_d}$ and
$\g(J)=g_{j_1}^{\b_1}\cdots g_{j_d}^{\b_d}$, where $\b_l=\b_l(I)$.
\begin{Lemma}
Let $n=3$. Assume $\b(f)=\g^{-1}\theta^{\m}\b(I)\g$, where 
$\m\in \Z, \g\in B_3$, $I=(i_1,\ldots,i_d)$ with  $i_1,\ldots,i_d \ge 3$. 
Then, we have 
\begin{enumerate}
\item $\calL(f)=  (-1)^{d+1}\displaystyle{\sum_{J\in S(I)}}(-1)^{|J|}\P_{\b(f)}(a_n^{\m}c(J)^{\g})$.
\item For any $J\in S(I)$,  the coefficient of $\P_{\b(f)}(a_n^{\m}c(J)^{\g})$ in $\calL(f)$ is non-zero.
\end{enumerate}
\end{Lemma}

{\it Proof}. (i)
For a partition $\calB$, let $K(\calB)$ be the set of
$l\in \Z_d$ with $(l,l+1)\in \calB$. 
For $J\in \Z^d(I)$, let $L(J)$ be the set of $l\in \Z_d$ with $(j_{l},j_{l+1})=(i_l,2)$.
Also, let $\calP_J$ be the set of $\calB\in \calP(d)$ with $K(\calB)\subset L(J)$. 
For $l\in \Z_d$, we have  $W_I((l))=\sum_{j=2}^{i_l}g_j^{\b_l}$ and 
 $W_I((l,l+1))=g_{i_l}^{\b_l}a_2^{\b_{l+1}}=-g_{i_l}^{\b_l}g_2^{\b_{l+1}}$, where we put $\b_{l+1}=e$ if $l=d$.
Therefore, for any partition $\calB$, we have  
$W_I(\calB)=(-1)^{\sharp K(\calB)}\displaystyle{\sum_{J: \calB\in  \calP_J}\g(J)}$,
and hence 
\begin{equation*}\sum_{\calB\in \calP(d)}W_I(\calB)
=\sum_{\calB\in \calP(d)}\sum_{J: \calB\in  \calP_J}(-1)^{\sharp K(\calB)}\g(J)
=\sum_{J\in \Z^d(I)}\e(J) \g(J),\end{equation*}
where $\e(J)=\sum_{\calB\in\calP_J}(-1)^{\sharp K(\calB)}$.
If $J\in \Z^d(I)-S(I)$, then 
$\e(J)$ is equal to some mutiple of the sum of $(-1)^{\sharp A}$ over the subsets $A$ of  $L(J)$. Since
$L(J)$ is  not empty, this sum is equal to zero.
If $J\in S(I)$, $\calP_J$ consists of  a single partition $\{(1),\ldots,(d)\}$, and hence 
$\e(J)=1$.
Therefore, Theorem 1 and the equality $\g(J)=(-1)^{d+|J|}c(J)$ imply the equality (i).
 
(ii) can be proved similarly as Lemma 6(iii). 
\endproof

\section{Proof of Proposition 2}
We first show that it is enough for the proof to consider the case of $\m=0$,
namely the case where $\b(\f)$ is conjugate to $\b(I)$.
The reason is given as follows:
Note that the period $m(\f)$ does not depend on the choice of an isotopy $\{\f_t\}$,
but the braid $\b(\f)$ and the integer  $\n(\f)$ depend on it.
To clarify the dependence on an isotopy,
denote them by $\b(\f,\{\f_t\})$ and $\n(\f;\{\f_t\})$ respectively.
Let $R_t:D\to D$ be the rotation of the disk with angle $2\pi t$.
Then, if we denote by $\{\f_t'\}$ the composition of the isotopies $\{\f_t\}$ and $\{ R_{-\m t}\}$,
then $\b(\f,\{\f_t'\})$ is equal to $\theta^{-\m}\b(\f,\{\f_t\})$, and hence it is conjugate to $\b(I)$. 
Therefore, if the proposition is proved in the case of $\m=0$,
then $\n(\f,\{\f_t'\})=\LCM\{d,n-2\}/(n-2)$, and hence 
$\n(\f,\{\f_t\})=m(\f)\m+\n(\f,\{\f_t'\})=m(\f)\m+\LCM\{d,n-2\}/(n-2)$.
Thus, we can assume $\b(\f)=\g^{-1}\b(I)\g$ for some $\g\in B_n$.

Let $\bar{d}=\LCM\{d,n-2\}$, $m=\bar{d}/d$,  and   $\n=\bar{d}/(n-2)$.
We  shall  prove that the coefficient of  $\P_{\b(\f^m)}(a_n^{\n})$ in $\calL(\f^m)$ is non-zero.
Then, since $m$ and $\n$ are relatively prime, the assertion of the proposition follows from  Lemma 1.
Let $p$ be an integer with $1\le p\le d-n+3$ and let $q=p+n-3$. Then, we have 
\begin{equation} a_1^{\b_p}a_{n-1}^{\b_q}=a_1^{\b_p}(a_{n-1}^{\b(i_q)})^{\b_{q+1}}=a_1^{\b_{p}}(a_na_1^{-1})^{\b_{q+1}}
=a_1^{\b_{p}}a_n(a_1^{\b_{q+1}})^{-1},
\end{equation}
where we put $\b_{d+1}=\b_1$.

Assume that $n\ge 4$ and $i_1,\ldots,i_d\ge 2$.
Define integers $\bar{i}_1,\ldots,\bar{i}_{\bar{d}}$ by $\bar{i}_{l}=i_{[l]}$,
where $[l]$ is the integer with $1\le [l] \le  d$ and $[l] \equiv l$ modulo $d$.
Let $\bar{I}=(\bar{i}_1,\ldots,\bar{i}_{\bar{d}})$.   
By Lemma 6(iii), it is enough for the proof to show that 
 $\P_{\b(\f^m)}(a_n^{\n})=\P_{\b(\f^m)}(\l_{\calB}(\calJ)^{\g})$
 for some $\calB \in \calP(\bar{d})$  and some $\calJ\in S_{\calB}(\bar{I})$.
We see by Lemma 2 that  this equality is equivalent to 
\begin{equation}\P_{\b(\bar{I})}(a_n^{\n})=\P_{\b(\bar{I})}(\l_{\calB}(\calJ)).\end{equation}
For $1\le r \le \n$, let $p_r=(r-1)(n-2)+1,q_r=r(n-2)$ and $B_r=[p_r,q_r]$. 
Note that all of these blocks have length $n-2$, and 
$\{B_1,\ldots,B_\n\}$ is a partition of $\Z_{\bar{d}}$.
Let \begin{equation*}
 \calB_r=\begin{cases} \{B_r\} & \text{if} \quad \bar{i}_{p_r}\ge 3,\\
\{(p_r),[p_{r+1},q_r]\} & \text{if}\quad  \bar{i}_{p_r}=2.\end{cases}
\end{equation*}
Then $\calB_r$ is a partition of $\{p_r,\ldots,q_r\}$, and if we put
$\calB={\calB_1}\cup\cdots\cup{\calB_\n}$, we have ${\calB}\in \calP'(\bar{d})$.
Let $S_{\calB_r}(\bar{I})$ be $S_{B_r}(\bar{I})$ if $\bar{i}_{p_r}\ge 3$, 
and be $S_{(p_r)}(\bar{I})\times S_{[p_r+1,q_r]}(\bar{I})$
 if $\bar{i}_{p_r}=2$.

For $r=1,\ldots,\n$, let $\zeta_r=a_1^{\b_{p_r}}a_{n-1}^{\b_{q_r}}$.
We shall show that there exists an element $J_r$ of $S_{\calB_r}(\bar{I})$ with $\l_{\calB_r}(J_r)=\zeta_r$.
In the case of $\bar{i}_{p_r}\ge 3$, let $J_r=1$. Then, 
  $\l_{\calB_r}(J_r)=\l_{\calB_r}(1)=c_1^{\b_{p_r}}a_{n-1}^{\b_{q_r}}=\zeta_r$.
In the case of $\bar{i}_{p_r}=2$, let  $J_r=(0,0)$. 
Then, $\l_{\calB_r}(J_r)=\l_{(p)}(0)\l_{[p_{r+1},q_r]}(0)=a_2^{\b_{p_r}}a_{n-2}^{\b_{q_r}}$.
Since 
\[(a_1^{-1}a_2)^{\b(2)}=\left[ (a_3a_1^{-1})(a_2a_1^{-1})^{-1}(a_3a_1^{-1})^{-1}\right] a_3a_1^{-1}
=a_3a_2^{-1},\] 
and hence 
\begin{equation*} 
(a_1^{-1}a_2)^{\b_{p_r}}=((a_1^{-1}a_2)^{\b(2)})^{\b_{{p_r}+1}}=
(a_3a_2^{-1})^{\b_{{p_r}+1}}=(a_{n-1}a_{n-2}^{-1})^{\b_{q_r}},
\end{equation*}
 we have
\[a_2^{\b_{p_r}}a_{n-2}^{\b_{q_r}}=(a_1a_1^{-1}a_2)^{\b_{p_r}}a_{n-2}^{\b_{q_r}}=
a_1^{\b_{p_r}}(a_1^{-1}a_2)^{\b_{p_r}}a_{n-2}^{\b_{q_r}}=\zeta_r.\]
Therefore, $\l_{\calB_r}(J_r)=\zeta_r$.

Let $\calJ=(J_1,\ldots,J_{\n})$. Then, 
 $\calJ\in S_{\calB}(\bar{I})$ and, since $\l_{\calB_r}(J_r)=\zeta_r$, we have 
$\l_{\calB}(\calJ)=\zeta_1\cdots \zeta_{\n}$.
Applying (21) to each pair   $p_r, q_r$, we have 
\[\l_{\calB}(\calJ)=\prod_{r=1}^{\n} a_1^{\b_{p_r}}a_n(a_1^{\b_{{q}_r+1}})^{-1}.\]
Since $p_1=1, q_r+1=p_{r+1}, q_{\n}+1=\bar{d}+1$,  this is equal to $a_1^{\b_1}a_n^{\n}a_{1}^{-1}$.
Therefore, by (5), $\l_{\calB}(\calJ)$ is Reidemeister equivalent to
 $(a_1^{-1})^{\b_1} (a_1^{\b_1}a_n^{\n})=a_n^{\n}$
with respect to  $\b(\bar{I})$. Therefore, (22) is proved.

Assume that  $n=3$ and $i_1,\ldots,i_d\ge 3$.
In this case, $\bar{d}=d=\n$ and $m=1$.
Using Lemma 7(ii) and Lemma 2,  we see that  
it is enough for the proof to show that $\P_{\b(I)}(a_n^{\n})=\P_{\b(I)}(c(J))$ for some $J\in S(I)$.
Let $J=(3,\ldots,3)\in S(I)$. Then, by (21),  
\[c(J)= c_3^{\b_1}\cdots c_3^{\b_{d}}
=\prod_{l=1}^d a_1^{\b_{l}}a_3(a_1^{\b_{l+1}})^{-1}
=a_1^{\b_1}a_3^{d}a_{1}^{-1}.\]
Therefore, $\P_{\b(I)}(c(J))=\P_{\b(I)}(a_1^{\b_1}a_3^{d}a_1^{-1})= \P_{\b(I)}(a_3^{d})$.
Since $d=\n$, the proof is completed.

\section{Bounds for the Nielsen number}

As a byproduct of  Lemma 6 and Lemma 7, 
we can obtain the following upper and lower bounds for the Nielsen number $N(f)$.
\begin{Theorem}
Assume that $\b(f)$ is conjugate to $\theta^{\m}\b(I)$.
\begin{enumerate}
\item  If $n\ge 4$ and  $i_1,\ldots,i_d \ge 2$, then
 \[\sum_{\calB\in\calP'(d)}\sharp S_{\calB}(I)-(2n-2)
\le N(f)\le \sum_{\calB\in\calP'(d)}\sharp S_{\calB}(I).\]  
\item If $n=3$ and  $i_1,\ldots,i_d \ge 3$, then
$\sharp S(I)-4\le N(f)\le \sharp S(I)$.
\end{enumerate}
\end{Theorem}

{\it Proof}.  We prove (i). 
 Let  $\Psi:\G(I)\to \calR(\b(f))$ be the map defined by 
$\Psi((\calB,\calJ))=\P_{\b(f)}(a_n^{\m}\l_{\calB}(\calJ)^{\g})$.
Let $\G'(I)$ be the set of $(\calB,\calJ)\in \G(I)$ with $\sharp [(\calB,\calJ)]>1$.
Let $\calR'(\b(f))$ be the set of Reidemeister classes
$\a$ with $\fix_\a(f)$ having index less than $-1$. We shall show that $\calR'(\b(f))$ coincides with the image of $\G'(I)$ under $\Psi$.
As we have shown in the proof of Lemma 6 (iii),
  the coefficient  $n(\calB,\calJ)$ of $\Psi((\calB,\calJ))$ in $\calL(f)$
is equal to $(-1)^{d+1+e(\l_{\calB}(\calJ))}\sharp  [(\calB,\calJ)]$.
On the other hand, $n(\calB,\calJ)$ is  equal to $\ind(\fix_{\Psi((\calB,\calJ))}(f))$
by its definition. Hence, we have 
\begin{equation} \ind(\fix_{\Psi((\calB,\calJ))}(f))=(-1)^{d+1+e(\l_{\calB}(\calJ))}\sharp [(\calB,\calJ)].
\end{equation}
This implies that $\Psi((\calB,\calJ))\in \calR'(\b(f))$ if and only if
$(-1)^{d+1+e(\l_{\calB}(\calJ))}=-1$ and $\sharp [(\calB,\calJ)]>1$.
The former condition $(-1)^{d+1+e(\l_{\calB}(\calJ))}=-1$ is redundant, 
since the index of any fixed point class of $f$ is less than two (Jiang and Guo \cite{jg}).
Thus we have proved the equality  $\calR'(\b(f))=\Psi(\G'(I))$.

 $\G'(I)$ is a disjoint union of equivalence classes $[(\calB_1,\calJ_1)],\ldots, [(\calB_m,\calJ_m)]$,
where  $m=\sharp \Psi(\G'(I))$. We have the following inequality due to \cite{jg}  
(see the proof of Theorem 4.1 there):
\[\sum_{\a\in \calR'(\b(f))}(\ind(\fix_\a(f))+1)\ge 2\chi(D_n)=2-2n.\] 
This inequality and (23) imply that 
$2-2n\le \sum_{i=1}^m  (-\sharp [(\calB_i,\calJ_i)]+1)=-\sharp \G'(I)+m$, and hence we have
 $\sharp \Psi(\G'(I))=m\ge \sharp \G'(I)+2-2n$. Therefore, since $\Psi$ is injective on $\G(I)-\G'(I)$, 
we have by Lemma 6(iii) that
\[N(f)=\sharp \Psi(\G(I))=\sharp (\G(I)-\G'(I))+\sharp \Psi(\G'(I))\ge 
\sharp \G(I)+2-2n.\]
Also, it is obvious that $N(f)\le \sharp \G(I)$.
Since $\sharp \G(I)=\sum_{\calB\in \calP'(d)}\sharp S_{\calB}(I)$, we have the bounds in (i).

(ii) can be proved similarly.
\endproof


\end{document}